\documentclass{article}
\usepackage[utf8]{inputenc}
\usepackage{enumerate}
\usepackage{amsthm,amsmath,amssymb}
\usepackage{graphicx}
\usepackage{lineno}
\usepackage[colorlinks=true,citecolor=black,linkcolor=black,urlcolor=blue]{hyperref}
\usepackage[english]{babel}
\usepackage{amsfonts}
\usepackage{epsfig, subfigure}
\usepackage{amscd,latexsym}

\usepackage{enumitem}

\theoremstyle{plain}
\newtheorem{theorem}{Theorem}
\newtheorem{lemma}[theorem]{Lemma}
\newtheorem{cor}[theorem]{Corollary}
\newtheorem{prop}[theorem]{Proposition}

{\itshape}{\rmfamily}
\newtheorem{remark}[theorem]{Remark}

\newtheorem{defi}[theorem]{Definition}

\newcommand{\edim}{{\rm edim}}
\newcommand{\sdim}{{\rm sdim}}
\newcommand{\ldim}{{\rm ldim}}
\newcommand{\mdim}{{\rm mdim}}
\newcommand{\ddim}{{\rm ddim}}
\newcommand{\dmd}{{{\rm dmd}}}

\newcommand{\ex}{{\rm \lambda}}

\newcommand{\bit}{\begin{itemize}}
\newcommand{\eit}{\end{itemize}}
\newcommand{\ben}{\begin{enumerate}}
\newcommand{\een}{\end{enumerate}}

\newcommand{\colred}{\color{black}}

\title{Metric Locations in Pseudotrees: A survey and new results}

\author{Jos\'e C\'aceres~\footnote{Departamento de Matem\'{a}ticas, Universidad de Almer\'{\i}a, Spain (\texttt{jcaceres@ual.es})} \and Ignacio M. Pelayo~\footnote{Departament de Matem\`{a}tiques, Universitat Polit\`{e}cnica de Catalunya, Spain (\texttt{ignacio.m.pelayo@upc.edu})}}

\begin{document}
\maketitle

\begin{abstract}
Since the publication of the original definition of metric location, many other versions with different features have emerged. Typically, the computation of a minimum location set is NP-hard, so it is sensible to study their behavior for different graph families. The aim of this paper is to revise the literature on different metric locations in the families of paths, cycles, trees and unicyclic graphs, as well as, providing several new results on that matter.
\end{abstract}





\vspace{+.1cm}\noindent \textbf{Keywords:} Metric dimension, metric locating sets, distance in graphs, unyciclic graphs.

\vspace{+.1cm}\noindent \textbf{AMS subject classification:} 05C05, 05C12, 05C35, 05C75.



\section{Introduction}\label{sec1:intro}

Location problems consist of determining a reference set in a graph such that every vertex is unequivocally associated to a set of ``coordinates'' that localize it. 
Since the vertex set of a graph is enough for this task, the question is not to prove its existence, but to find the minimum one. 

Metric location in graphs were officially come into being with the publication of two seminal papers by P. Slater~\cite{s75} and F. Harary and R. A. Melter~\cite{hm76}. 
However, the field can be traced back to a paper by R. Silverman~\cite{s60} about the problem of detecting a certain subset of a given set by comparing it with some other subsets. 
The approach to solve that problem was to define a distance between the different subsets of a set and to choose carefully the comparisons needed. 
To the best of our knowledge, in this work the concept of metric basis is named for the first time.

The adjective ``metric'' is used through the paper {\colred to refer to those locations} in which the usual distance between graph vertices plays a crucial role, in the spirit of the original definition given in~\cite{s75,hm76}. By contrast, there exist other types of locations in which the role of the graph distance is substituted by the relation of adjacency of vertices, given rise to ``neighbor'' locations which are not treated here.

Since the publication of the works by P. Slater~\cite{s75} and F. Harary and R. A. Melter~\cite{hm76}, metric location theory has evolved at great velocity giving rise to a hundred of different variations, developing many new methods and deepening into several directions, although the first of those variations had to wait 20 years to appear in a paper published in 2003 by R. C. Brighham, G. Chartrand, R. D. Dutton and P. Zhang~\cite{bcdz03} who imposed the condition of being dominating to the location sets. 
Almost immediately, in 2004, A. Seb\H{o} and E. Tannier~\cite{st04} published their definition of strong metric dimension by using a certain notion of betweenness in geodetic paths. Just a year later, J.  Cáceres, C. Hernando, M. Mora, I. M. Pelayo, M. L. Puertas, C. Seara and D. R. Wood~\cite{chmppsw07} defined what is  currently  known as the  doubly metric dimension in which vertices have not only different coordinates but also, the coordinates of two vertices have not a fixed difference. 
C. Hernando, M. Mora, P. Slater and D. R.  Wood~\cite{hmsw08} asked themselves in 2008 what if an element of the metric location set ``fails" or give incorrect lectures, and gave rise to the fault-tolerant metric dimension. In 2010, F. Okamoto, B. Phinezy and P. Zhang {\colred introduced} in~\cite{opz10} a local metric dimension to distinguish only adjacent vertices instead of using a coloring. In 2015, A. Estrada-Moreno, I.G. Yero and J.A. Rodríguez-Velázquez in~\cite{eyr14} generalized the fault-tolerant dimension into the k-metric dimension by considering the case in which more than one element of the location set might fail. In 2017, A. Kelenc, D. Kuziak, A. Taranenko and I.G. Yero~\cite{kkty17} distinguished vertices and edges with the mixed metric dimension, and a year later, A. Kelenc, N. Tratnik and I.G. Yero~\cite{kty18} introduced the edge metric dimension for locating only edges.

One of the crucial facts of the field is that computing a minimum metric location set is typically NP-hard. 
Therefore, it becomes essential to identify graph families for which this computation can be accomplished efficiently in polynomial time. 
Among these families, trees are often considered, and various algorithms have been developed specifically for this purpose.

Since trees do not contain cycles, it is natural to explore the graph family characterized by a single cycle: unicyclic graphs. 
There exists an extensive volume of research papers that makes it exceedingly difficult to follow the advancements in this particular area. 
Consequently, the authors recognized the necessity of conducting a comprehensive survey on the topic of metric location within this graph family. 
Although there are many location parameters, including variations of variations, we opted to choose nine of them having in common that the corresponding dimensions for paths, cycles and trees are well-known (see Table~\ref{megatable1}). 
However, not in all cases a complete characterization of locations in unicyclic graphs has been found and we were able to provide some new results in four of the dimensions studied. 

The paper is organized as follows: after concluding this section by a description of the notation that will be followed in the rest of the work, Section~\ref{dmd} is dedicated to doubly metric locating sets where the work has been completed. 
This metric location parameter is chosen to be the first  because of the application it has for the rest of metric locations. 
Section~\ref{dim} is devoted to the original and { \colred first metric location} that gave rise to the entire field. 
Here, we provided significant new results as well as for the  strong metric dimension which is studied in Section~\ref{sdim}. 
We finish with the dominating metric dimension in Section~\ref{ddim} where again the goal has been achieved.

In two of the remaining five parameters (mixed and local metric dimensions) other authors have proved bounds and exact values (see Table~\ref{megatable2}).
However, for the edge, fault-tolerant and k-metric dimensions (see Table~\ref{megatable3}) there exist partial results, but not a complete characterization of those dimensions for unicyclic graphs. 
In Section~\ref{omlp}, we offer a glimpse of the actual landscape for those cases and the paper is finished with a section of conclusions and further work.

\subsection{Basic terminology}

All the graphs considered are undirected, simple, finite and (unless otherwise stated) connected.
Unless otherwise specified,  $G=(V,E)$ stands for a graph of order $n$ and size {\colred $m$}, being $|V|=n$, $V=[n]=\{1,\ldots,n\}$ and $|E|=m$.
Let $v$ be a vertex of $G$.
The \emph{open neighborhood} of $v$ is $\displaystyle N_G(v)=\{w \in V(G) :vw \in E\}$, and the \emph{closed neighborhood} of $v$ is $N_G[v]=N_G(v)\cup \{v\}$ (we will write $N(v)$ and $N[v]$ if the graph $G$ is clear from the context).
The \emph{degree} of $v$ is $\deg(v)=|N(v)|$.
The minimum degree  (resp., maximum degree) of $G$ is $\delta(G)=\min\{\deg(u):u \in V(G)\}$ (resp., $\Delta(G)=\max\{\deg(u):u \in V(G)\}$).

For any two vertices $u,v\in V(G)$ of a connected graph $G$, a $u-v$ \emph{geodesic} is a  $u-v$ shortest  path, i.e., a $u-v$ path of minimum order. 
The length of a $u-v$ geodesic is called the \emph{distance} $d_G(u,v)$ between $u$ and $v$, or simply  $d(u,v)$, when the graph $G$ is clear from the context. 
The diameter of $G$ is ${\rm diam}(G) = \max\{d(v,w) : v,w \in V(G)\}$.

Let $S=\{w_1,w_2,\ldots,w_k\}$ a set of vertices of a graph $G$.
The distance between a vertex $v\in V(G)$ and  $S$, denoted by  $d(v,S)$, is the minimum of the distances between $v$ and the vertices of $S$, that is, $d(v,S)=\min\{d(v,w):w\in S\}$.
Given a vertex $v \in V(G)$ and an edge $e = xy \in E(G)$, the distance\ between $v$ and $e$ is 
$d_G(v, e) = d_G(v,\{x,y\})$.
The \emph{metric representation} $r(v|S)$ of a vertex $v$ with respect to $S$ is defined as the $k$-vector $r(v|S)=(d_{G}(v,w_1),d_{G}(v,w_2)\ldots,d_{G}(v,w_k))$.

Let $u,v \in V(G)$ be  a pair of vertices such that  $d(u,w)=d(v,w)$ for all $w\in V(G)\setminus\{u,v\}$, i.e.,  such that  either $N(u)=N(v)$ or $N[u]=N[v]$.
In both cases, $u$ and $v$ are said to be \emph{twins}.
Let $W$ be a set of vertices of $G$.
If the vertices of $W$ are pairwise twins, then $W$ is called a \emph{twin set} of $G$.

Let $W\subseteq V(G)$ be a subset of vertices of  $G$.
The  \emph{closed neighborhood} of $W$ is $N[W]=\cup_{v\in W} N[v]$.
The subgraph of $G$ induced by $W$, denoted by $G[W]$, has $W$ as vertex set and $E(G[W]) = \{vw \in E(G) : v \in W,w \in W\}$.

\begin{figure}[ht]
\begin{center}
\includegraphics[width=0.35\textwidth]{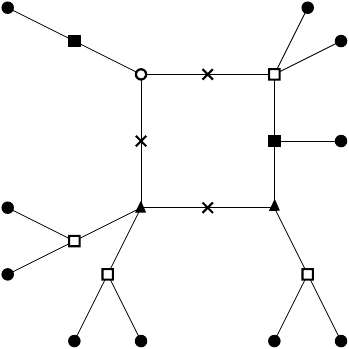}
\caption{ (White) square vertices: (Strong) support vertices.
The white circle is an exterior major vertex but it is neither strong nor a support vertex.
Triangle vertices are both major vertices but none of them is neither an exterior major vertex nor a support vertex.
Check that $g=8$, $\ell=10$, $\lambda=6$, $\ell_s=8$, $\lambda_s=4$, $c_2=3$, $c_3=5$, $\rho=3$.}
\label{niceunic}
\end{center}
\end{figure}

For a graph $G$ of minimum degree $\delta(G)=1$, end-vertices are called \emph{leaves}.
The set and the number of leaves of $G$ are denoted by ${\cal L}(G)$ and $\ell(G)$, respectively.
A \emph{support vertex} (resp., \emph{strong support vertex}) is a vertex adjacent to a leaf (resp., to at least two leaves).
A \emph{major vertex} of $G$ is any vertex of degree at least 3.
A \emph{terminal vertex} of a major vertex $v$ of $G$ is a leaf $u$ such that 
$d(u,v)<d(u,w)$ for every other major vertex $w$ of $G$.
An \emph{exterior major vertex}  
(resp., \emph{strong exterior major vertex})
is a major vertex  with at least one terminal vertex (resp., at least two terminal vertices)
(see Figure \ref{niceunic}).

The number of exterior major vertices (resp., strong exterior major vertices) of $G$ is denoted by $\lambda(G)$ (resp., $\lambda_s(G)$).
A leaf of $G$ is called \emph{strong} if its exterior major vertex is strong.
The number of strong leaves of $G$ is denoted  by $\ell_s(G)$.
Clearly, $\ell_s(G)-\lambda_s(G)=\ell(G)-\lambda(G)$.

A graph $G$ of order $n$ and size $m$ is called \emph{unicyclic} if $n=m$.
{\colred
In the sequel (unless stated otherwise), if $g$ is the girth of $G$, then its unique cycle will be denoted by $C_g$.
}
A unicyclic graph $G$ of order $n$ and girth $g$ is said to be \emph{proper} if $n>g$, i.e. if it contains at least a leaf.
If $g$ is odd (resp., even), then $G$ is said to be an \emph{odd (resp., even) unicyclic graph}.
The connected component of $G-E(C_g)$ containing a vertex $v\in V(C_g)$  is denoted by $T_v$ and  is called the \emph{branching tree} of $v$.
The tree $T_v$ is said to be \emph{trivial} if $V(T_v)=\{v\}$.
A branching tree is called a \emph{thread} if both $T_v$ is a path and $\deg(v)=3$.
A vertex $v\in V(G)$ is a \emph{branching vertex} if  either $v \not\in V (C_g)$ and $\deg (v) \ge 3$ or $v \in V (C_g)$ and $\deg (v ) \ge 4$. 
A vertex $v \in V (C_g)$ is \emph{branch-active} if $T_v$ contains a branching vertex.
The number of branch-active  vertices  of $C_g$ is denoted by $\rho(G)$.

A triple of {\colred different} vertices $u,v,w \in V(C_g)$ is called a \emph{geodesic triple} of $G$ if
$d(u,v)+d(v,w)+d(w,u)=g$.
A vertex $v$ from $C_g$ is called a \emph{root vertex} on $C_g$ if $T_v$ is non-trivial and otherwise,  $v$ is said to be a \emph{trivial vertex} of $C_g$. 
The set and the number of all root vertices of $C_g$  is denoted by $C_3(C_g)$ and $c_3(C_g)$, while the set and the number of all trivial vertices of $C_g$  are denoted by $C_2(C_g)$ and $c_2(C_g)$, respectively. 
Certainly, every vertex of $C_g$ of  degree at least 3 is a root vertex, while all vertices of $C_g$ that are of degree 2 are trivial.

%
%


A set $U$ of vertices of a graph $G$ is called \emph{independent} if no two vertices in $U$ are adjacent. 
The \emph{independence number} of $G$, denoted by $\alpha(G)$, is the cardinality of a maximum independent set of $G$.

A vertex $v$ of a graph $G$ is said to be a \emph{boundary vertex} 
of a vertex $u$ if no neighbor of $v$ is further away from $u$ than $v$, i.e., if for every vertex $w\in N(v)$, 
$d(u,w) \le d(u,v)$. 
The set of boundary vertices of a vertex $u$ is denoted by $\partial(u)$.
The \emph{boundary} of $G$, denoted by $\partial(G)$, is the set of all of its boundary vertices, i.e., 
$\partial(G)=\cup_{u \in V(G)}\partial(u)$.
Given a  pair of vertices $u,v\in V(G)$ if $v\in \partial(u)$, then $v$ is also said to be \emph{maximally distant} from $u$.
Moreover, a pair of  vertices $u,v\in V(G)$ are said to be \emph{mutually maximally distant}, or simply $MMD$, if both $v\in \partial(u)$ and $u\in \partial(v)$.
Notice that, as was pointed out in \cite{ryko14}, the boundary of $G$ can also be defined as the set of $MMD$ vertices of $G$, i.e., 
$\partial(G)=\{v\in V(G): {\rm there \, exists} \,  u\in V(G) \, {\rm such \, that} \, u,v \,{\, \rm are} \,  MMD \}$.
A pair of $MMD$  vertices of a cycle are also called \emph{antipodal}.

For additional details and information on basic graph theory, we refer the reader
to \cite{clz16}.

\section{Doubly metric dimension}\label{dmd}

This parameter was formally introduced by J. C\'aceres et al. in \cite{chmppsw07} as a means to study the metric dimension in grids, cylinders and tori. 

In this paper, it will also be  a powerful tool to study unicyclic graphs, so we commence with this section.  

Let $G=(V,E)$ be a graph.
A pair of vertices $u, v\in V$ is said to \emph{doubly resolve} two vertices $x,y \in V$ if
$$d_G(x,u)-d_G(x,v) \neq d_G(y,u)-d_G(y,v).$$

\begin{defi} \label{dls} 
{\rm A  set of vertices $S$ of  $G$ is called \emph{doubly locating} if every pair of {\colred distinct} vertices $x,y \in V$ is doubly resolved by some pair of vertices of $S$.}
\end{defi}

A doubly locating set of  minimum cardinality is called  a \emph{doubly metric basis} of $G$. 
The \emph{doubly metric dimension} of $G$, denoted by $\dmd(G)$,  is the  cardinality of a doubly metric basis. 
To know more about this parameter see mainly  \cite{chmppsw07,j22.2} and also \cite{kkcs12.1,kkcs12.2,mkkc12}.

The problem of computing the doubly metric dimension of trees was approached and completely solved in \cite{chmppsw07}.

\begin{theorem}
{\rm\cite{chmppsw07}}
Let  $T$ be  a tree  having $\ell(T)$ leaves. 
Then, $\dmd(T)  = \ell(T)$.
\label{dmd.trees}
\end{theorem}

In particular, $\dmd(P_n)  = 2$.
The problem of characterizing all minimum doubly locating sets of a cycle $C_g$ was  implicitly approached and  partially proved in  \cite{chmppsw07}.
{\colred
Here we present the whole result.
}

\begin{theorem}
Consider the cycle $C_g$ and assume that $V(C_g)=[g]$.

\begin{enumerate}[label=\rm \bf(\arabic*)]

\item If $g$ is odd, then $\dmd(C_g)=2$.
\newline
Moreover, a pair of vertices is a doubly metric basis if and only if they are antipodal.

\item
If $g$ is even, then $\dmd(C_g)=3$.
\newline
Moreover, a triple of vertices is a doubly metric basis if and only if it is a geodesic triple.

\end{enumerate}
\label{dmd.cycles}
\end{theorem}
\begin{proof}

\begin{enumerate}[label=\rm \bf(\arabic*)]

\item
Suppose that $g=2k+1$ and take the antipodal pair $\{1,k+1\}$.
Notice that for every $i \in [g]$:

$$ d(i,k+1)-d(i,1)=   
\left \{ \begin{array}{cl}
{\colred k-2i+2} & {\rm if}\, 1 \le i \le k+1, \\
2i-3k-3 & {\rm if}\, k+2 \le i \le 2k+1.  \\
\end{array}\right . $$

Let $\{i,j\}$ such that $1 \le i < j \le 2k+1$.
We distinguish cases:

{\bf Case (a)}:
If $1 \le i < j \le k+1$, then 
{\colred 
$$d(i,k+1)-d(i,1) = k-2i+2 \neq k-2j+2 = d(j,k+1)-d(j,1).$$
}

\vspace{-.7cm}
{\bf Case (b)}:
If $k+2 \le i < j \le 2k+1$, then 
$$d(i,k+1)-d(i,1) = 2i-3k-3 \neq 2j-3k-3 = d(j,k+1)-d(j,1).$$

\vspace{-.3cm}
{\bf Case (c)}:
If $1 \le i \le k+1$ and $k+2 \le  j \le 2k+1$, then 
{\colred 
$$d(i,k+1)-d(i,1) = k-2i+2 \neq 2j-3k-3 = d(j,k+1)-d(j,1),$$ 
\vspace{-.06cm}
since otherwise, if $k-2i = 2j-3k-3$,then $2i+2j = 4k+5 = 2i+2j \ge 4k+4$, a contradiction.
}

To end the proof of this item, consider a pair of non-antipodal vertices $\{1,h\}$.
W.l.o.g., we can assume that $1 < h < k+1$.
Notice that the pair $\{h,h+1\}$ is not doubly resolved by $\{1,h\}$, since $d(h,h)-d(h,1)=1-h=d(h+1,h)-d(h+1,1)$.

\item
Suppose that $g=2k$ and take a pair of vertices $\{1,h\}$.
W.l.o.g., we can assume that $1 < h \le k$.
If $h<k$, then the pair $\{k,k+1\}$ is not doubly resolved by $\{1,h\}$, 
since $d(k,h)-d(k,1)=1-h=d(k+1,h)-d(k+1,1)$.
If $h=k$, then the pair $\{k-1,k+1\}$ is not doubly resolved by $\{1,k\}$, 
since $d(k-1,k)-d(k-1,1)=3-k=d(k+1,k)-d(k+1,1)$.
Hence, $\dmd(C_g)\ge3$.

Take a geodesic triple of vertices $S$.
W.l.o.g., we can assume that $S=\{1,\alpha, \beta\}$ and $1 < \alpha < k \le \beta  \le 2k-1$.
Take a pair pf vertices $\{i,j\}$.
W.l.o.g., we can assume that $1 < i < \alpha$.
We distinguish cases:

{\bf Case (a)}:
If $1 < i < j < \alpha$, then $\{1,\alpha\}$ doubly resolves $\{i,j\}$, since
$d(i,\alpha)-d(i,1) = \alpha-2i+1 \neq \alpha-2j+1 = d(j,\alpha)-d(j,1)$.

{\bf Case (b)}:
If $\alpha < j \le k$, then 
$\{1,\alpha\}$ doubly resolves $\{i,j\}$, since
$d(i,\alpha)-d(i,1) = \alpha-2i+1 \neq -\alpha +1 = d(j,\alpha)-d(j,1)$.

{\bf Case (c)}:
If $k < j$, then 
$\{1,\beta\}$ doubly resolves $\{i,j\}$.
To prove this we distinguish cases.

{\bf Case (c.1)}:
If $i+k-1 \le \beta $ and $j < \beta$, then
$d(i,\beta)-d(i,1)=(i-1+2k-\beta)-(i-1)=2k-\beta$ and $d(j,\beta)-d(j,1)=(\beta-j)-(2k-j)=\beta-2k$.

{\bf Case (c.2)}:
If $i+k-1 \le \beta $ and $\beta < j$, then
$d(i,\beta)-d(i,1)=(i-1+2k-\beta)-(i-1)=2k-\beta$ and $d(j,\beta)-d(j,1)=(j-\beta)-(2k-j)=2j-\beta-2k$.
Notice that if  $2k-\beta=2j-\beta-2k$, then $j=2k$, a contradiction.

{\bf Case (c.3)}:
If $\beta < i+k-1$ and $j<\beta$, then
$d(i,\beta)-d(i,1)=(\beta-i)-(i-1)=\beta-2i+1$ and $d(j,\beta)-d(j,1)=(\beta-j)-(2k-j)=\beta-2k$.
Notice that if $\beta-2i+1=\beta-2k$, then $2i-1=2k$, a contradiction.

{\bf Case (c.4)}:
If $\beta < i+k-1$ and $ \beta < j$, then
$d(i,\beta)-d(i,1)=(\beta-i)-(i-1)=\beta-2i+1$ and $d(j,\beta)-d(j,1)=(j-\beta)-(2k-j)=2j-\beta-2k$.
Notice that if $\beta - 2i + 1 = 2j-\beta - 2k$, then $ 2i + 2j = 2\beta + 2k + 1$, a contradiction.

Hence, we have proved both that $\dmd(C_g)=3$ and that every geodesic triple is a doubly metric basis.
To end the proof of this item, consider a non-geodesic triple of vertices $S=\{ 1, \alpha, \beta \}$.
W.l.o.g., we can assume that $1 < \alpha < \beta < k$.
Notice that the pair $\{\beta,\beta+1\}$ is not doubly resolved neither by $\{1,\alpha\}$, nor by $\{1,\beta\}$, nor by $\{\alpha, \beta\}$, since for every $x\in S$, $d(\beta+1,x)=d(\beta,x)+1$.
\end{enumerate}
\end{proof}

In \cite{j22.2}, after noticing that every doubly locating set of a unicyclic graph $G$ must contain all its leaves, it was proved that if 
{ \colred
$G$ is a proper uniclyclic graph of girth $g$,  then  
$\ell(G) \le \dmd(G) \le \ell(G)+1$ if $g$ is odd, and $\ell(G) \le \dmd(G) \le \ell(G)+2$ if $g$ is even.
}
Starting both from this result and from Theorem \ref{dmd.cycles}, we have characterized the different families of unicyclic  graphs achieving each of the possible values for its doubly metric dimension.

{\colred
\begin{lemma}
Let $u,v$ a pair of  vertices of the cycle $C_g$ of a unicyclic graph $G$.
Let $\ell_u$ (resp., $\ell_v$) either a leaf of the branching tree $T_u$  (resp., the branching tree  $T_v$) or the vertex $u$ (resp., $v$). 
If $\{x,y\} \subsetneq V(C_g)$ is doubly resolved by $\{u,v\}$, then it is also doubly resolved by $\{\ell_u,\ell_v\}$.
\label{lem:dmd.unic1}
\end{lemma}
\begin{proof}
Suppose that $d(u,\ell_u)=k_u$ and $d(u,\ell_u)=k_v$.
If $d(x,\ell_u)-d(x,\ell_v)=h_u$ and $d(y,u)-d(y,v)=h_v$, then
$d(x,\ell_u)-d(x,\ell_v)=h_u+k_u-k_v \neq h_v+k_u-k_v =d(y,\ell_u)-d(y,\ell_v)$, since $h_u \neq h_v$.
\end{proof}
}

\begin{theorem}
{\colred Let  $G$ be  a proper unicyclic graph such that $\ell(G)=\ell$. }
Then, 

\begin{enumerate}[label=\rm \bf(\arabic*)]

\item If $g$ is odd, then 
$ \dmd(G)  =  \left \{ \begin{array}{cl}
\ell  & {\rm if}\, C_3(C_g)\, {\rm contains \, an\, antipodal\, pair}, \\
\ell + 1 & {\rm otherwise}. \\  
\end{array}\right . $

\item
If $g$ is even, then
$ \dmd(G)  =  \left \{ \begin{array}{cl}
\ell  & {\rm if}\, C_3(C_g)\, {\rm contains \, a\, geodesic\, triple}, \\
\ell + 2 & {\rm if}\, c_3(C_g)=1, \\
\ell + 1 & {\rm otherwise}. \\  
\end{array}\right . $

\end{enumerate}
\label{dmd.unic1}
\end{theorem}
\begin{proof}

\begin{enumerate}[label=\rm \bf(\arabic*)]

\item
If $C_3(C_g)$ contains an antipodal pair then, {\colred according to Lemma \ref{lem:dmd.unic1} and } Theorem \ref{dmd.cycles} {\bf(1)}, the set ${\cal L}(G)$ of leaves of $G$ is a doubly  locating set of $G$.

Conversely, if $C_3(C_g)$ contains no antipodal pair, then to obtain a doubly  locating set of $G$, in addition to ${\cal L}(G)$, a trivial vertex $h$  of $C_g$ must be added, to be more precise, a vertex $h$ of $C_g$ antipodal to one of the vertices of $C_3(C_g)$.
Then,  {\colred according to Lemma \ref{lem:dmd.unic1} and }   Theorem \ref{dmd.cycles} {\bf(1)}, the set ${\cal L}(G)\cup \{h\}$  is a doubly locating set of $G$.

\item
If $C_3(C_g)$ contains a geodesic triple then, {\colred according to Lemma \ref{lem:dmd.unic1} and }  Theorem \ref{dmd.cycles} {\bf(2)}, the set ${\cal L}(G)$ of leaves of $G$ is a doubly  locating set of $G$.

Conversely, suppose that $C_3(C_g)$ contains no geodesic triple. 
If $c_3(C_g)=1$ and $C_3(C_g)=\{h_1\}$, take a pair of trivial vertices $h_2,h_3 \in [g]$ such that that $\{h_1,h_2,h_3\}$ be a geodesic triple of $C_g$.
Then, {\colred according to Lemma \ref{lem:dmd.unic1} and }  Theorem \ref{dmd.cycles} {\bf(2)}, the set ${\cal L}(G)\cup \{h_2,h_3\}$  is a minimum doubly locating set of $G$.
If $c_3(C_g)\ge 2$ and $\{h_1, h_2\}\subseteq C_3(C_g)$, take a trivial vertex $h_3 \in [g]$ such that that $\{h_1,h_2,h_3\}$ be a geodesic triple of $C_g$.
Then, {\colred according to Lemma \ref{lem:dmd.unic1} and }  Theorem \ref{dmd.cycles} {\bf(2)}, the set ${\cal L}(G)\cup \{h_3\}$  is a minimum doubly locating set of $G$.

\end{enumerate}
\end{proof}

Hence, in this section we have completed previous works to obtain a total characterization of the doubly metric dimension of unicyclic graphs.

\section{Metric dimension}\label{dim}

Being introduced primarily in~\cite{hm76,s75}, a large part of the literature regarding with graph locating have to do with metric dimension, the first and more important way of locating vertices in a graph. 

Let $G=(V,E)$ be a graph.
A vertex $v\in V$ is said to \emph{resolve} two vertices $x,y \in V$ if $d_G(x,v)\neq d_G(y,v)$.

\begin{defi} \label{ls} 
{\rm A  set of vertices $S$ of  $G$ is called \emph{metric-locating} if  every pair of distinct  vertices  $x,y \in V$ is resolved by a vertex of $S$.}
\end{defi}

Metric-locating sets are also known as locating sets and resolving sets.
Moreover, if $S$ is a metric-locating set of a graph $G$, then it is usually said that $S$ resolves the set of vertices $V(G)$.

A metric-locating set of  minimum cardinality is called  a \emph{metric basis} of $G$. 
The \emph{metric dimension} of $G$, denoted by $\dim(G)$,  is the  cardinality of a metric basis. 
To know more about this parameter see mainly  \cite{chmppsw07, cejo00, hm76,s75} and also \cite{bcggmmp13, bc11, bm11, bdfhmp18, bcpz03,  cgh08, eky17,hmpscp05, hmpsw10, krr96, pz02, st04, ss21.5,ss22,ss22.2,ss22.1}.

It is straightforward to check, first,  that $\dim(P_n)=1$, being a vertex $v$  of $P_n$ a metric basis if and only if $v$ is a leaf, and second, that  $\dim(C_g)=2$, being any pair of vertices (resp., of non-antipodal vertices) of $C_g$ a metric basis if $g=2k+1$ is odd (resp., $g=2k$ is even).

\begin{lemma}\cite{cejo00}
\label{lem:cejo00}
If $G$ is a graph, then $\dim(G)\geq \ell(G)-\ex(G)$.
\end{lemma}

Moreover, it was also proved in \cite{cejo00} that any locating set of a graph $G$ must contain, at least, all but one of the terminal vertices of  every  exterior major vertex of $G$.

\begin{theorem}
{\rm\cite{cejo00,hm76,s75}}
{\colred
For every tree $T$, $\dim(T)  = \ell(T)-\lambda(T)$.
}
\label{thmc3.f1}
\end{theorem}

At this point, we are interested in finding the metric dimension of unicyclic graphs. In~\cite{ss21.5,ss22}, the authors gave a characterization for {\colred the metric dimension of  unicyclic graphs  based} on forbidden configurations. Here, we present a different approach that relies on the parameters of the graph and the use of doubly resolving sets. Naturally, our results could be derived from those presented in~\cite{ss21.5,ss22} and vice versa. However, in our opinion, this parameter approach has some value, and that providing different lines of attack to the same problem would enrich the discussion. 

\begin{theorem}
{\rm \cite{ss21.5}}
{\colred
Let  $G$  be a unicyclic graph  such that $\ell(G)=\ell$, $\lambda(G)= \lambda$  and $\rho(G)=\rho$.}
If $\hat{\rho}=max\{2-\rho,0\}$, then  
$$\ell - \lambda + \hat{\rho} \le \dim(G ) \le  \ell - \lambda + \hat{\rho} +1.$$
\label{cryptic}
\end{theorem}
In the previous section, it was talked about the use of doubly locating sets for obtaining results in metric dimension. In the case of unicyclic graphs, the relationship is given by the next result.

\begin{lemma}
\label{lem:dr4md}
Let  $G$ be  a proper  unicyclic graph of girth $g$. 
If $S$ is a doubly locating set of $V(C_g)$, then $S$ resolves $V(C_g)$ together with the vertices of all threads of $G$.
\end{lemma}
\begin{proof}
Let $v,v'\in S$. 
Assume that there exists one thread beginning in a certain vertex $u_i\in V(C_g)$. Let $\{x_0,x_1,...,x_r\}$ the vertices of the thread in consecutive order where $x_0=u_i$. 
Clearly, $r(x_k|S)=r(u_i|S)+(k,k,\ldots,k)$ and so, two vertices of the same thread cannot have the same metric representation and thus, they are resolved by $S$. 

However, it is possible that $r(u_j|S)=r(x_k|S)$ for $u_j\in V(C_g)$ and $x_k$ in the thread. 
In this case $r(u_j|S)=r(x_k|S)=r(u_i|S)+(k,k,\ldots,k)$ and that means that $d(v,u_j)-d(v,u_i)=d(v',u_j)-d(v',u_i)=k$ contradicting the fact that $S$ is a doubly locating set for the vertices of the cycle. 

Finally, suppose that there are at least two threads $\{x_0,x_1,...,x_r\}$ and $\{y_0,y_1,...,y_s\}$ where $u_i=x_0$ and $u_j=y_0$ and such that $r(x_k|S)=r(y_l|S)$ (w.l.o.g we can assume that $k\geq l$). 
Thus, 
{\colred 
$r(u_i|S)+(k,k,\ldots,k)=r(y_l|S)=r(u_j|S)+(l,l,\ldots,l)$
}
and this implies that 
$d(v,u_i)-d(v,u_j)=d(v',u_i)-d(v',u_j)=k-l$ and again $S$ would not doubly resolve $u_i$ and $u_j$. 
So, $S$ resolves the vertices of the cycle together with the vertices of any thread.
\end{proof}

In the rest of the section, we divide the study into two parts depending on the parity of $C_g$. 

\subsection{Odd unicyclic graphs}\label{dim.odduni}

Let us start with the easiest case, although a complete explicit characterization seems to be  very complicated and obscure. 
The cases  in which $\rho(G)=0,1$ have been solved with the following two results.

\begin{lemma} 
Let  $G$ be  an odd proper unicyclic graph.
If $\rho(G)=0$, then $\dim(G)=2$.
\label{lem:rho=0.odd}
\end{lemma}
\begin{proof}
Let $v$ and $v'$ be a pair of antipodal vertices of $C_g$. 
By Theorem~\ref{dmd.cycles}, they doubly resolve all the vertices of the cycle $C_g$, and hence, according to Lemma~\ref{lem:dr4md}, since $\rho(G)=0$, they  resolve all the vertices of $G$.
\end{proof}

\begin{lemma} 
Let  $G$ be  an odd proper unicyclic graph. 
If $\rho(G)=1$, then $\dim(G)=\ell(G)-\lambda(G)+1$.
\label{lem:rho=1.odd}
\end{lemma}
\begin{proof}
Under the condition of the hypothesis, there is exactly  one branch-active vertex $v$ in the cycle $C_g$, being $T_v$  its branching tree. 

Given a metric basis $S$ of $G$ and according to Lemma~\ref{lem:cejo00}, all the terminal vertices except one of each exterior major vertex with positive terminal degree of $T_v$ should belong to $S$. 
However, since $v$ is a cut-point of $G$, not all the vertices of $S$ lie in $T_v$ otherwise $S$ will not resolve two vertices of $C_g$ at the same distance from $v$. 
Thus, a lower bound for a metric dimension is $\dim(G)\geq \ell(G)-\lambda(G)+1$.

Now, let $S$ be a set with $\ell(G)-\lambda(G)$ leaves together with a vertex $v'$ antipodal of $v$. If $u\in S$ is a leaf, then it is clear that $\{u,v'\}$ doubly resolves all the vertices of the cycle, and thus by Lemma~\ref{lem:dr4md}, $S$ resolves all the vertices of the cycle and the threads. 

Moreover, $S$ resolves all the vertices in $T_v$ except perhaps $v$ and other vertex $u'$ with the same exterior major vertex as $v$. 
But in this case, they are resolved by $v'$ since $d(v,v')<d(u',v')$. 
Finally, a vertex in $T_v$ and any other vertex in $V(G)\setminus V(T_v)$ is resolved by a leaf in $S$.

Thus, $S$ is a metric-locating set achieving the lower bound $\ell(G)-\lambda(G)+1$, i.e., $S$  is a metric basis of $G$.
\end{proof}

For the rest of the cases, i.e., when $\rho(G)\geq 2$,  we can give a partial result of the value of the metric dimension.

\begin{lemma}
Let  $G$ be  an odd proper unicyclic graph of girth $g$. 
If $g\ge 3$, $\rho(G)\ge 2$ and there are, at least, two antipodal branch-active vertices, 
then $\dim(G)=\ell(G)-\lambda(G)$. Particularly, if $g=3$ and $\rho(G)\geq 2$ then $\dim(G)=\ell(G)-\lambda(G)$.
\label{lem:rho>=2.odd}
\end{lemma}
\begin{proof}
Let $S$ be a vertex subset formed by $\ell-\lambda$ leaves for each branching tree. 
By Lemma~\ref{lem:cejo00},  $\dim(G)\geq |S|$, and thus it remains to prove that $S$ is a metric-locating set. 
We first claim  that $S$ doubly resolves the vertices of $V(C_g)$.

Let $T_1,T_2$ be the branching trees whose branch-active vertices $u_1,u_2$ are antipodal, and let $v_1,v_2$ be two leaves such that $v_i\in V(T_i)\cap S$ for { \colred $i\in\{1,2\}$}. 
By Lemma~\ref{dmd.cycles}, $\{u_1,u_2\}$ doubly resolves the vertices of cycle. 
{ \colred Moreover}, for any $x\in V(C_g)$, we have that $r(x|\{v_1,v_2\})=r(x|\{u_1,u_2\})+(d(u_1,v_1),d(u_2,v_2))$, and hence $\{v_1,v_2\}$ also doubly resolves the vertices of the cycle, and so do $S$ since $\{v_1,v_2\}\subseteq S$.

By Lemma~\ref{lem:dr4md}, $S$ resolves all the vertices of the cycle and the threads of $G$. 
A similar reasoning can be used to prove that $S$ distinguishes between a vertex in a branching-tree, and a vertex either in a thread or in the cycle $C_g$.

Two vertices in the same branching-tree $T$ different from its branch-active vertex are resolved by the leaves in $V(T)\cap S$. 
A vertex in $T$ and its branch-active vertex $u$ are resolved by a leaf of other branching-tree in $S$. So $S$ distinguishes two vertices of the same branching-tree.

Finally, let us consider two vertices $x$ and $y$ in different branching-trees. 
Its braching-active vertices are joined by a shortest path in the cycle. 
If we delete the rest of the cycle and all the vertices hanging out, we build a subtree $T$ of $G$ which is also an isometric subgraph. 
By Theorem~\ref{thmc3.f1}, $x$ and $y$ are resolved in $T$ by the vertices in $V(T)\cap S$. 
However, since $T$ is an isometric subgraph, $x$ and $y$ are also resolved in $G$ by the same vertices in $S$. And this completes the proof.

For the particular case of $g=3$, note that if $\rho(G)\geq 2$ there are always two branching-active vertices that are antipodal.
\end{proof}

{\colred
The metric dimension of an odd unicyclic graph $G$ with at least two 
branch-active vertices without antipodal branch-active vertices may be, 
according to Theorem \ref{cryptic}, either  $\dim(G)=\ell(G)- \lambda(G) $ or $\dim(G)=\ell(G) - \lambda(G) +1 $. 
}

\begin{figure}[htbp]
\begin{center}
\includegraphics[width=10cm]{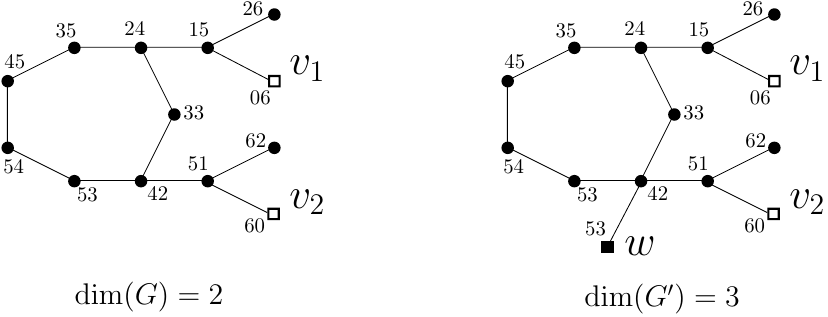}
\end{center}
\caption{Left: Odd unicyclic graph $G$ with $\ell(G)=4$, $\lambda(G)=2$, $\rho(G) =2$, without antipodal branch-active vertices. 
Right: Odd unicyclic graph $G'$ with  $\ell(G)=5$, $\lambda(G)=3$, $\rho(G) =2$,  without antipodal branch-active vertices.}
\label{fig:odd.rho>=2}
\end{figure}

As an example, the unicyclic graph $G$ displayed in  Figure~\ref{fig:odd.rho>=2} (left) has $\dim(G)= \ell(G)-\lambda(G)$ (a metric basis is $\{v_1,v_2\}$). Nevertheless, if  $G'$ is the graph obtained from $G$ by adding  a vertex $w$ as shown in Figure~\ref{fig:odd.rho>=2} (right), then  $\dim(G')=\ell(G')-\lambda(G')+1$, being  a metric basis  the set $\{v_1,v_2,w\}$.

As an immediate consequence of the previous 3 lemmas, the following result is obtained.

\begin{cor}\label{cor:xulo.odd}
{\colred
Let  $G$ be  an odd unicyclic graph such that $\ell(G)=\ell$, $\lambda(G)= \lambda$  and $\rho(G)=\rho$.
}
If $g\ge 3$, then $\dim(G)=\ell - \lambda + \max\{2-\rho,0\} $ whenever any the following conditions holds.

\begin{enumerate}[label=\rm \bf(\arabic*)]

\item $0 \le \rho \le 1$.

\item $\rho\ge 2$ and there are, at least, two antipodal branch-active vertices.

\item $g=3$.

\end{enumerate}
\end{cor}

\subsection{Even unicyclic graphs}\label{dim.evenuni}

The case in which the unicyclic graph has an even cycle turns out to be much more difficult than the previous case. 
Hence, the results are not as complete as in the previous case.

\begin{figure}[htbp]
\begin{center}
\includegraphics[width=12cm]{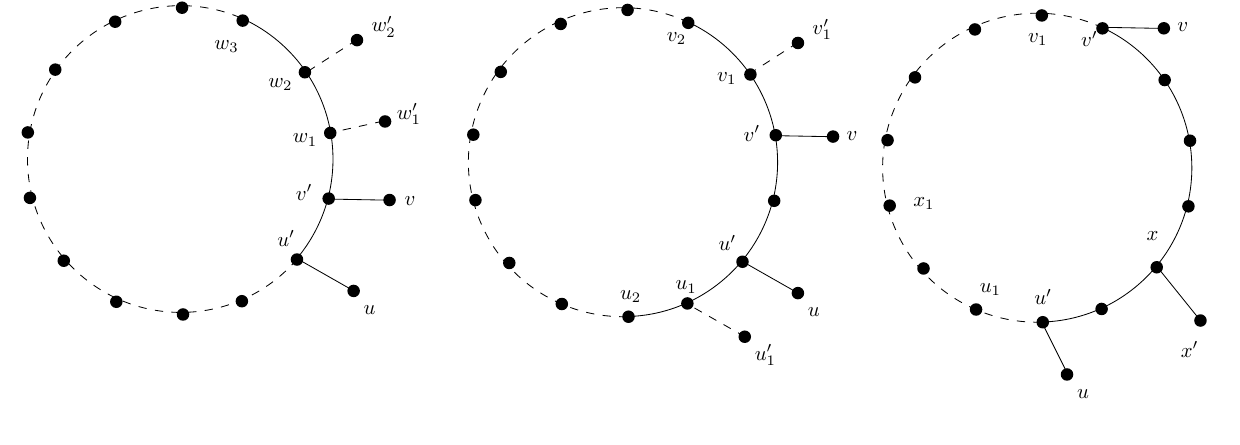}
\end{center}
\caption{Several situations of the proof of Proposition~\ref{lem:evencycle.rho=0}}
\label{fig:cycleeven}
\end{figure}

\begin{prop}\label{lem:evencycle.rho=0}
{\colred
Let  $G$ be  an even proper unicyclic graph.
}
If $g\ge8$, $\rho(G)=0$ and $c_2(C_g)\le1$, then $\dim(G)=3$.
\end{prop}
\begin{proof}
On the contrary, suppose that $\{u,v\}$ is a metric basis of $G$ in the conditions of the hypothesis. 
Since $u$, $v$ or both may be a vertex in a thread or in the cycle, we will use the following notation for treating all the cases at once: 
If $u$ is a vertex in the interior of a thread, thn $u'$ will be the vertex of that thread that also belongs to the cycle, otherwise $u=u'$. 
Analogously, we  define $v'$.

Note that $u'$ and $v'$ cannot be antipodal in the cycle since otherwise $\{u,v\}$ will not be resolving. 
Thus, there is a unique shortest path between $u'$ and $v'$ and $d(u',v')\leq k-1$. 
We  organize the proof by studying different distances between $u'$ and $v'$. 

Firstly, suppose that $d(u',v')=1$ (see Figure~\ref{fig:cycleeven} left). 
Since $g\geq 8$, we can add $w_1,w_2,w_3\in V(C_g)$ to obtain the shortest path $u-v-w_1-w_2-w_3$. 
The minimum number of threads is {\colred $2k-1=g-1$}, so $w_1$ or $w_2$ is the beginning of a thread and if we denote $w_i'$ {\colred($i\in\{1,2\}$)} as the adjacent to $w_i$ in the thread then we have that $r(w_i'|\{u,v\})=r(w_{i+1}|\{u,v\})$ and hence $\{u,v\}$ is not a locating set. 

Consider now the case $2\leq d(u',v')<k-1$ (see Figure~\ref{fig:cycleeven} center). 
Again since $g\geq 8$, we can take $u_1,u_2,v_1,v_2\in V(C_g)$ to construct the path $u_2-u_1-u'-\ldots -v'-v_1-v_2$. 
This path might not be a shortest path but the subpaths $u'-v_2$ and $u_2-v'$ are shortest paths and that is enough for our reasoning. 
Then, either $u_1$ or $v_1$ is the beginning of a thread. 
Without loss of generality, we can assume that it is $v_1$ and its adjacent in the thread is $v_1'$. 
Then, $r(v_1'|\{u,v\})=r(v_2|\{u,v\})$ and again $\{u,v\}$ is not a locating set.

The final case occurs when $d(u',v')=k-1$ (see Figure~\ref{fig:cycleeven} right). Let $u_1$ and $v_1$ be the adjacent vertices to $u'$ and $v'$ respectively in the cycle which are not in the shortest path $u'-v'$ (see Figure~\ref{fig:cycleeven}). Then the shortest path $u_1-v_1$ has the same length as $u'-v'$, and for every interior vertex $x$ of $u'-v'$, corresponds with $x_1$ in $u_1-v_1$ such that $d(u',x)=d(u_1,x_1)$. 
Moreover, $r(x_1|\{u,v\})=r(x|\{u,v\})+(1,1)$. 
Since there are $k-2$ internal vertices in the path $u'-v'$, at least one of them, say $x$, has a thread and its adjacent vertex in that thread, say $x'$, has the same metric representation with respect to $\{u,v\}$ as $x_1$. 
So again $\{u,v\}$ cannot be a locating set and $\dim(G)=3$.
\end{proof}

\begin{figure}[htbp]
\begin{center}
\includegraphics[width=5cm]{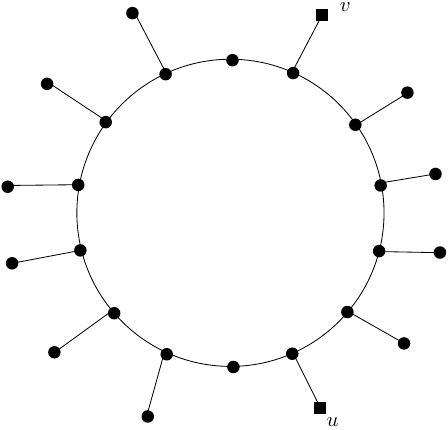}
\end{center}
\caption{Even unicyclic graph $G$ with $g-2$ threads and $\dim(G)=2$.}
\label{fig:example}
\end{figure}

The above result is the best bound we can obtain, since there are unicyclic graphs such that the number of threads is $2k-2$ and the metric dimension equals two. 
Consider for example a unicyclic graph $G$ with girth $g=2k$ with $k\geq 3$ and such that every vertex has a thread of length one except two antipodal vertices. 
In Figure~\ref{fig:example}, it is shown the case $k=7$,  where the pair $\{u,v\}$ of square vertices is a metric basis. 
The construction can be repeated for any integer $k\geq 3$.

Now let us check {\colred some} cases we left out in the previous result.

\begin{prop}
{\colred
Let $G$ be an uniciclic graph such that $\rho(G)=0$.
}
\begin{enumerate}[label=\rm \bf(\arabic*)]
\item If $g=4$, then $\dim(G)=2$.
\item If $g=6$, then $c_2(C_g)\neq 0$ if and only if $\dim(G)=2$.
\end{enumerate}
\end{prop}
\begin{proof}
\begin{enumerate}[label=\rm \bf(\arabic*)]
\item Consider first a unicyclic graph $G$ with $g=4$ and $\rho(G)=0$. 
Whenever you have two threads whose initial vertices are adjacent, say $u_1,u_2$, in the cycle then pick the leaves of those two threads as metric basis, say $v_1,v_2$. 
If the only two threads of $G$ have as initial vertices $u_1$ and $u_3$, two non-adjacent vertices of the cycle, then pick $\{v_1,u_2\}$ as a metric basis. 

If $G$ contains only one thread in, say $u_1$, then $\{v_1,u_2\}$ will be a metric basis. Finally, if $G$ does not have threads then $G=C_4$ and two consecutive vertices as $u_1,u_2$ are enough to resolve the graph.

\item Suppose now that $g=6$. 
If $c_2(C_g)=0$, then there exists a leaf for every vertex in the cycle. 
The reader can easily check that 
{\colred
the metric dimension of the corona product of $C_g$ and $K_1$ is  $\dim(C_g\odot K_1)=3$
} 
and so $\dim(G)=3$. 
If $c_2(C_g)\neq 0$, then there exists at least one vertex in the cycle without a thread. 
If there is only one $u_1$, then pick the leaves hanging from the adjacent vertices in the cycle to $u_1$. 
It is not difficult to check that those pair of vertices form a metric basis. 
The same applies if there is a vertex in the cycle adjacent to two vertices with threads. 

Otherwise, we have at  least two adjacent vertices $u_1,u_2\in V(C_g)$ such that other adjacent to one of them, say $u_3$, has a thread. Let $v_3$ be the leaf in that thread. Again it is not difficult to prove that $\{u_1,v_3\}$ is a basis for $G$ and therefore $\dim(G)=2$.
\end{enumerate}
\end{proof}

The next result is very similar to the Lemma~\ref{lem:rho>=2.odd}, with the particularities of dealing with an even cycle.

\begin{prop}\label{lem:evencycle.rho>=3}
{\colred
Let  $G$ be  an even unicyclic graph.
}
If $\rho(G)\geq 3$ and there are three branching-active vertices forming a geodesic triple, then 
$\dim(G)=\ell(G)-\lambda(G)$.
\end{prop}
\begin{proof}
Let $S$ be a vertex subset with $\ell(G)-\lambda(G)$ leaves. 
By Lemma~\ref{lem:cejo00}, we know that any locating set must contain $S$, and thus $\dim(G)\geq |S|$. 
We next show that $S$ is a locating set.

First, we claim that a subset of $S$, and consequently $S$, doubly resolves the vertices of the cycle. Let $T_1,T_2,T_3$ be the trees whose branching-active vertices $u_1,u_2,u_3$ is a geodesic triple. We construct $S'=\{v_1,v_2,v_3\}$ where $v_i\in V(T_i)\cap S$, so $S'\subseteq S$. As $\{u_1,u_2,u_3\}$ is a geodesic triple, then they doubly resolves the vertices of the cycle, by Theorem~\ref{dmd.cycles}. That means:
\begin{multline}
\forall x,y\in V(C_g), \exists i,j\in \{1,2,3\}: d(x,u_i)-d(x,u_j)\neq d(y,u_i)-d(y,u_j)\Rightarrow\\
\Rightarrow d(x,u_i)-d(x,u_j)+(d(u_i,v_i)-d(u_j,v_j))\neq d(y,u_i)-d(y,u_j)+(d(u_i,v_i)-d(u_j,v_j))\Rightarrow\\
\Rightarrow d(x,u_i)+d(u_i,v_i)-(d(x,u_j)+d(u_j,v_j))\neq d(y,u_i)+d(u_i,v_i)-(d(y,u_j)+d(u_j,v_j))\Rightarrow\\
\Rightarrow d(x,v_i)-d(x,v_j)\neq d(y,v_i)-d(y,v_j)
\end{multline}
And consequently, $S'$ also doubly resolves the vertices of the cycle. 
The rest of the proof is similar to the one of Lemma~\ref{lem:rho>=2.odd}.
\end{proof}

\begin{figure}[htbp]
\begin{center}
\includegraphics[width=10cm]{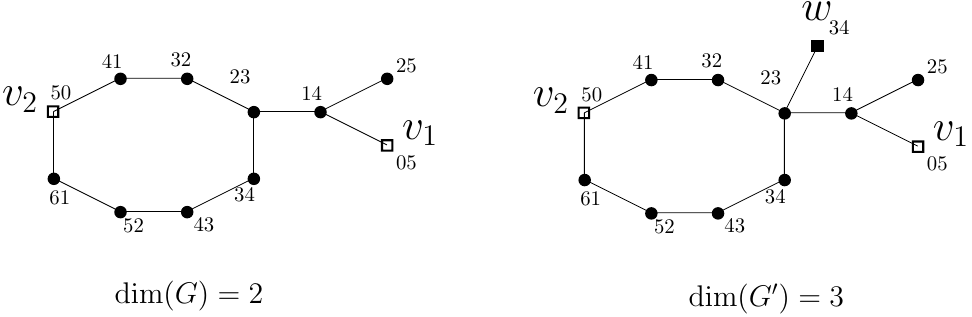}
\end{center}
\caption{Left: Even unicyclic graph $G$ with $\ell(G)=2$, $\lambda(G)=1$, $\rho(G) =1$. 
Right: Even unicyclic graph $G'$ with  $\ell(G')=3$, $\lambda(G')=2$, $\rho(G') =1$.}
\label{fig:even.rho=1}
\end{figure}

As an immediate consequence of the previous 3 lemmas, the following result is obtained.

\begin{cor}\label{cor:xulo.even}
{\colred
Let  $G$ be  an even unicyclic graph such that $\ell()=\ell$, $\lambda(G)=\lambda$, $\rho(G)=\rho$  and $c_2(C_g)=c_2$.
}
If $g\ge 3$, then $\dim(G)=\ell - \lambda + \max\{2-\rho,0\} $ whenever any the following conditions holds.

\begin{enumerate}[label=\rm \bf(\arabic*)]
\item $\rho=0, g\geq 8$ and $c_2\geq 2$.
\item $\rho=0$ and $g=4$.
\item $\rho=0, g=6$ and $c_2\neq 0$.
\item $\rho\geq 3$ with three branching-active vertices forming a geodesic triple.
\end{enumerate}
\end{cor}

\section{Strong metric dimension}\label{sdim}
The strong metric dimension is the most atypical one among the metric dimensions. To start with, it does not rely in coordinates to distinguish the vertices of the graph but in a certain notion of betweenness in shortest paths.

Let $G=(V,E)$ be a graph. 
A vertex $v\in V$ is said to \emph{strong resolve} two vertices $x,y \in V(G)$ if there exists either  a $v-x$ geodesic that contains $y$, or a $v-y$ geodesic  that contains $x$. 

\begin{defi} \label{sls} 
{\rm A  set of vertices $S$ of $G$ is called \emph{strong locating} if every pair of distinct vertices $x,y \in V(G)$ is strong resolved by a vertex $v \in S$.}
\end{defi}

In other words, $S$ is a strong locating set if, for every pair of distinct vertices $x,y \in V(G)$, there is a vertex $v \in S$, such that either 
$d_G(v,x)=d_G(v,y)+d_G(y,x)$ or $d_G(v,y)=d_G(v,x)+d_G(x,y)$.

A strong locating set of  minimum cardinality is called  a \emph{strong metric basis} of $G$. The \emph{strong metric dimension} of $G$, denoted by $\sdim(G)$,  is the  cardinality of a strong metric basis. 
This parameter was formally introduced  by A. S\H{e}bo and E. Tannier in \cite{st04}. 
To know more about this parameter see mainly  \cite{k20,op07,ryko14,st04} and also \cite{gsw14,kkcm14,kkcs12.1,kkcs12.2,kyr13,mo11,y13}.

The next result is mentioned in many papers, and mainly by this reason, we include here the proof for the sake of completeness.

\begin{theorem}
Let $S$  be a proper subset of vertices  of a graph $G=(V,E)$.
Then, the following statements are equivalent.

\begin{enumerate}[label=\rm \bf(\arabic*)]

\item $S$ is a strong locating set og $G$.

\item $G$ is uniquely determined by the matrix of distances between the vertices of $S$ and the vertices of $V$.


\end{enumerate}
\label{sdim.dmatrix}
\end{theorem}
\begin{proof}
Assume that $S$ is a strong locating set of $G$.
Take $u,v \in V \setminus S$.
Let $w \in S$ such that $d(w,v)=d(w,u)+d(u,v)$.
We distinguish cases.

\vspace{.2cm}
{\bf Case (a)}:
If $d(w,v) \ge d(w,u)+2$, then $d(u,v) \ge 2$, i.e., $uv \not \in E)$.

\vspace{.2cm}
{\bf Case (b)}:
If $d(w,v) = d(w,u)+1$, then $d(u,v) = 1$, i.e., $uv  \in E)$.

\vspace{.2cm}
Conversely, suppose that $S$ is not a strong locating set of $G$.
Take $u,v \in V \setminus S$ such that for every vertex $w \in S$, 
$d(w,v) < d(w,u) + d(u,v)$ and $d(w,u) < d(w,v) + d(v,u)$.

Take $w \in S$ and $z \in  V \setminus S$.
If  $\eta$ is a $w-z$ geodesic, then either $u$ or $v$ does not belong to $\eta$, since otherwise either
$d(w,v) = d(w,u) + d(u,v)$ or  $d(w,u) = d(w,v) + d(v,u)$.
Hence, both if $uv  \in E$ and $uv \not \in E$, the matrix of distances between the vertices of $S$ and the vertices of $V$ is the same.
\end{proof}

As it was mentioned above, strong locating sets do not rely on coordinates so it would be desirable to find a characterization of some of the vertices that belong to a strong locating set. 
The next pair of results gives us precisely that characterization.

\begin{lemma}\label{sdimlema1}
Let $S$ be a strong locating set of a graph $G$. 
If $u,v \in V(G)$ are \emph{MMD}, then $\{u,v\}\cap S \neq \emptyset$.
\end{lemma}
\begin{proof}
On the contrary, assume $u,v\notin S$.
Let $w\in S$ such that the pair $u,v$ is strongly resolved by $w$. 
Suppose that $u$ lies  in a $w-v$ geodesic. 
Thus, there exists a neighbor of $u$ in this path that is farther away from $v$ than $u$, contradicting the fact that $u,v$ are \emph{MMD}.
\end{proof}

\begin{theorem}
Consider the cycle $C_g$.
Then, $\sdim(C_g)=\lceil \frac{g}{2} \rceil$.
Moreover,
a $\lceil \frac{g}{2} \rceil$-set $S$ is a  strong metric basis of $C_g$ if and only if
the set $V(G)\setminus S$ contains no  antipodal pair.
\end{theorem}
\begin{proof}
Assume that $V(C_g)=[g]$ and take the set $W=[\lceil \frac{g}{2} \rceil]$.
Notice that for every pair of vertices $i,j$ such that $\lceil \frac{g}{2} \rceil < i < j \le g$, 
$ d(i,1) = d(i,j) + d(j,1)$, which means that there is a (unique) $(i-1)$ geodesic containing  vertex $j$.
Hence, $W$ is a strong locating set of $C_g$ of cardinality $\lceil \frac{g}{2} \rceil$.

Let $S$ be a strong locating set of $C_g$.
According to Lemma~\ref{sdimlema1}, the set $V(C_g)\setminus S$ contains no  antipodal pair, which means that the cardinality of $V(G)\setminus S$ is at most $\lfloor \frac{g}{2} \rfloor$.
Hence, $\sdim(C_g) = g -\lfloor \frac{g}{2} \rfloor = \lceil \frac{g}{2} \rceil$.

Conversely, let $S$ be a $\lceil \frac{g}{2} \rceil$-set such that $V(G)\setminus S$ contains no  antipodal pair.
Let  $i,j\in V(C_g)\setminus S$ such that $d(i,j)=h < \lfloor \frac{g}{2} \rfloor$.
W.l.o.g., we may assume that $i=1$ and $j=h+1$.
If $\lfloor \frac{g}{2} \rfloor=k$, take the vertex $k+1$, which is an antipodal vertex of $1$, and check that 
the there is a (unique) $1-(k+1)$ geodesic containing vertex $h$.
\end{proof}

\begin{lemma}
Let $U$ be a $k$-set of non-consecutive vertices of the odd cycle $C_{2k+1}$.
Then, it contains at least a pair of antipodal vertices.
\end{lemma}
\begin{proof}
Suppose that $U$ contains no antipodal pair.
Let $\{U_i\}^{\alpha}_{i=1}$ be the maximal disjoint partition of $U$ such that, for every $i\in[\alpha]$, $U_i$  is a $\beta_i$-set consisting of $\beta_i$ consecutive vertices.
If $W_i$ denotes the set of antipodal vertices of $U_i$, then $\{W_i\}^{\alpha}_{i=1}$ is a disjoint partition contained in $V(C_g)\setminus U$ such that, for every $i\in[\alpha]$, $W_i$  is a $(\beta_i+1)$-set consisting of $\beta_i+1$ consecutive vertices.
Then, $\sum_{i=1}^{\alpha}\beta_i=k$, $\Big| \bigcup_{i=1}^{\alpha} W_i\Big| \le k+1$,
$\Big|\bigcup_{i=1}^{\alpha} W_i \Big| = \sum_{i=1}^{\alpha}(\beta_i+1) =\sum_{i=1}^{\alpha}\beta_i +\alpha = k +\alpha$, and thus $\alpha=1$, i.e., $U$ is  a $k$-set of consecutive vertices.
\end{proof}

As a straight consequence of both  the previous Theorem and this Lemma, the following result holds.

\begin{cor}
Every strong metric basis of the cycle odd  $C_{2k+1}$ consists of $k+1$ consecutive vertices.
\end{cor}

\begin{theorem}{\rm \cite{st04}}~
Let $T$ be a tree with $\ell$ leaves.
Then, any set consisting of all the leaves except one is a minimum strong locating set of $T$.
Hence, $\sdim(T)=\ell-1$.
\end{theorem}
\label{sdim.trees}
\begin{proof}
According to Lemma~\ref{sdimlema1}, for any pair of vertices mutually maximally distant, at least one of them must belong to a given  strong locating code $S$. 
In a tree, the vertices which are  mutually maximally distant are the leaves, so $S$ must contain all the leaves except one. 

Conversely, let $S$ be a set containing all the leaves of $T$ except one. 
It is clear that the shortest path between any pair of vertices $u,v$  can be extended until both of its extremes are leaves. 
Hence, at least one of them is in $S$ and strongly resolves the pair $u,v$.
\end{proof}

\begin{figure}[!htbp]
\begin{center}
\includegraphics[width=0.85\textwidth]{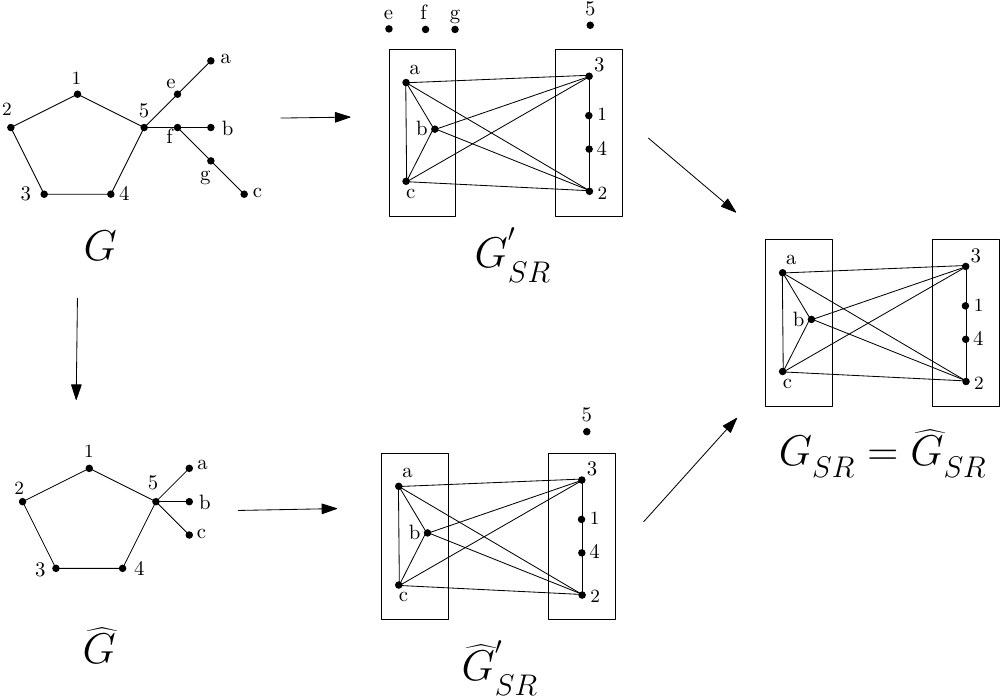}
\end{center}
\caption{A unicyclic graph $G$, its closed necklace $\widehat{G}$ and its strong resolving graph $G_{SR}$.}
\label{fig:sdim3}
\end{figure}

To compute the strong metric dimension, it is necessary, in most cases, to use the so-called  strong resolving graph introduced by O. R. Oellermann and J. Peters-Fransen in~\cite{op07}.
The \emph{strong resolving graph} $G^{'}_{SR}$ of a connected graph $G$ is the graph with vertex set 
{\colred $V(G'_{SR})=V(G)$} where two vertices are adjacent if they are mutually maximally distant in $G$.
Another definition, given in~\cite{ryko14}, reduces the vertex set to be $V(G_{SR})=\partial(G)$.
{\colred Thus, the strong resolving graph $G_{SR}$ defined in \cite{ryko14} is the subgraph of the strong resolving graph $G^{'}_{SR}$ defined in \cite{op07}  obtained from the latter graph by deleting its isolated vertices} (see Figure \ref{fig:sdim3}). 
Although we will stick formally with the second definition, this does not affect to the proof of any of the results given.

\begin{theorem}{\rm \cite{op07,ryko14}}
Let $G$ be a connected graph. 
Then,
$$ \sdim(G)=|\partial(G)|-\alpha(G_{SR}).$$
\label{sdim.partalpha}
\end{theorem}

Observe that, $[P_n]_{SR}=K_2$, $[C_n]_{SR}=C_n$ and $[T_{\ell}]_{SR}=K_{\ell}$, where $T_{\ell}$ denotes any tree with $\ell$ leaves.
Notice also that $\alpha(K_r)=1$ and  $\alpha(C_n)=\left \lfloor \frac{g}{2} \right \rfloor$.
Hence,  as a direct consequence of the previous theorem, it is derived that $\sdim(P_n)=1$, 
$\sdim(C_n)=\left \lceil \frac{n}{2} \right \rceil$ and $\sdim(T_{\ell})=\ell-1$.

The so-called \emph{closed necklace} $\widehat{G}$  associated to a unicyclic graph $G$ of girth $g$, is the unicyclic graph where, for every vertex $i \in V(C_g)$, the branching tree $T_i$ of $i$ is replaced by the 2-path $K_2$ if $T_i$ is a path, and by the star $K_{1,\ell_i}$, where $\ell_i$ is the number of leaves of $T_i$, otherwise.
It is a rutinary exercise to check, first, that $\partial(\widehat{G})=\partial(G)$, and second, that 
$\widehat{G}_{SR}=G_{SR}$ (see Figure \ref{fig:sdim3} for an example).

Let $G=(V,E)$ be a closed necklace of order $n$, girth $g$ with $\ell$ leaves.
Clearly, if $C_g$ is this cycle and $c_2(C_g)=h$, then $|\partial(G)|=\ell+h$.
Moreover, the vertex set $V(G_{SR})=\partial(G)$ of its strong resolving graph $G_{SR}$,
can be partitioned into two sets $U$ and $W$, 
in such a way  that $G_{SR}[U]=K_{\ell}$, and $G_{SR}[W]$ is a graph of order $h$ and maximum degree at most 2.

\begin{figure}[!htbp]
\begin{center}
\includegraphics[width=0.75\textwidth]{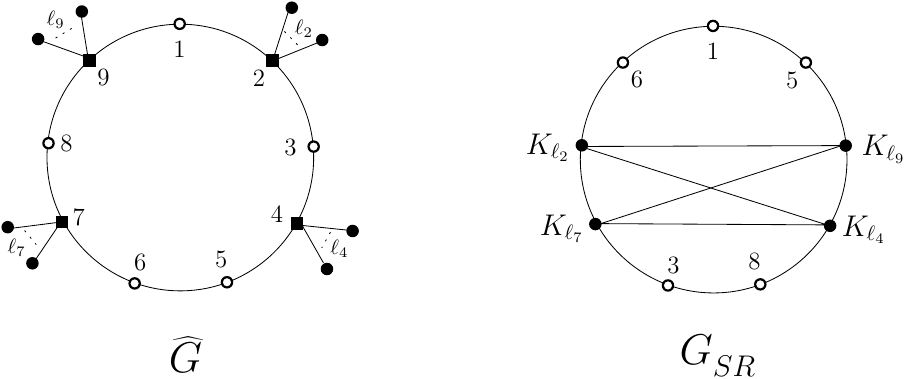}
\end{center}
\caption{Closed necklace $\widehat{G}$ of a unicyclic graph $G$ of girth 9 such that $\ell=\ell_2+\ell_4+\ell_7+\ell_9$,  $c_2(C_g)=5$, $G_{SR}[W]=P_3+P_2$ and its strong resolving graph $G_{SR}$.}
\label{fig:sdim5}
\end{figure}

\begin{remark}
If $g$ is odd, then  $G_{SR}$ is a Hamiltonian graph  and $G_{SR}[W]$ is a disjoint union of paths (see Figure \ref{fig:sdim5}).
Hence, $\displaystyle \left \lceil \frac{h}{2} \right \rceil \le \alpha(G_{SR}) \le h+1$.
\label{sdimremark1}
\end{remark}

\begin{figure}[!htbp]
\begin{center}
\includegraphics[width=0.72\textwidth]{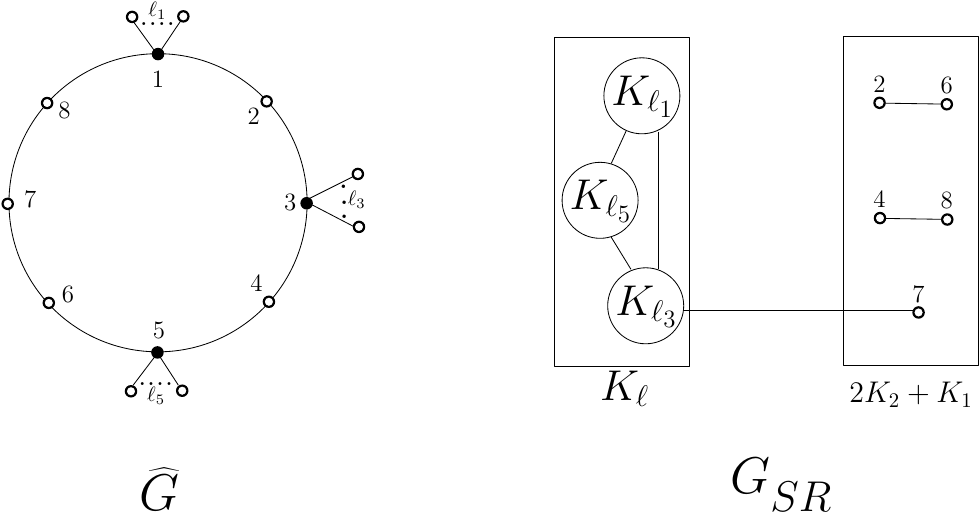}
\end{center}
\caption{Closed necklace $\widehat{G}$ of a unicyclic graph $G$ of girth 8 with $\ell=\ell_1+\ell_3+\ell_5$ leaves, $r=2$, $s=1$ and $t=1$ and its  strong resolving graph $G_{SR}$.}
\label{fig:sdim4}
\end{figure}

\begin{remark}
If  $g$ is even, then $G_{SR}[W]=rK_2+sK_1$, where $r$ is the number of antipodal pairs of $C_2(C_g)$, $t$ is the number of  antipodal pairs of $C_3(C_g)$ and $2s=g-2r-2t$ 
(see Figure \ref{fig:sdim4}).
Hence, $h=2r+s$ and 
$\displaystyle \left \lceil \frac{h}{2} \right \rceil = r + \left  \lceil \frac{s}{2} \right \rceil  \le r+s  \le  \alpha(G_{SR}) \le r+s+1 $.
\label{sdimremark2}
\end{remark}

The next result, proved by D. Kuziak in \cite{k20}, provides the range of values for the strong metric dimension of unicyclic graphs.

\begin{theorem}{\rm \cite{k20}}\label{thm:k20}
{\colred
Let  $G$ be  a proper  unicyclic graph of girth $g$ with $\ell$ leaves. 
}
If $c_2(C_g)=h$,  then 
$$
\max \left \{ \left \lceil \frac{g}{2} \right \rceil,\ell-1\right \}
\leq 
\sdim(G)
\leq 
\ell + \left \lfloor \frac{h}{2} \right \rfloor.$$
\end{theorem}
\label{sdim.bounds}
%

Now, we address the problem not only of proving that the  bounds displayed in Theorem~\ref{thm:k20} are tight but also of characterizing the families of graphs achieving each of them.

\begin{prop}\label{prop.st1}
{\colred
Let $G$ be a proper unicyclic graph of girth $g$ with $\ell$ leaves. 
}
If $c_2(C_g)=h$,  then 
$\sdim(G)=\ell-1$ if and only if either $g=3$ and $h=0$ or $g\ge4$ and any of the following conditions hold:

\begin{enumerate}[label=\rm \bf(\arabic*)]

\item $ 0 \le h \le 1$,

\item $g=2k$ even, $\displaystyle  2 \le h \le k-1$ and $C_2(C_g)$ contains no pair of antipodal vertices.

\item $g=2k+1$ odd, $\displaystyle  2 \le h \le k-1$, $C_2(C_g)$ contains no pair of antipodal vertices and $C_3(C_g)$ contains at least one antipodal triple.

\end{enumerate}
\label{sdim.ell-1}
\end{prop}
\begin{proof} 
According to Theorem  \ref{sdim.partalpha}, 
$\sdim(G) = |\partial(G)|-\alpha(G_{SR}) = \ell + h - \alpha(G_{SR})$.
Hence, $\sdim(G)=\ell-1$ if  and only if $\alpha(G_{SR})=h+1$.

If $g=3$ then, in all cases, $\alpha(G_{SR})=1$.
Hence, $\sdim(G)=\ell+h-1=\ell-1$ if and only if $h=0$.

Suppose that $g\ge4$.
We distinguish cases:

{\bf Case (a)}:
If $h=0$, then $G_{SR}\cong K_{\ell}$.
Thus, $\alpha(G_{SR})=1=h+1$.

{\bf Case (b)}:
If $h=1$, then $G_{SR}$  is a graph of order $\ell+1$ containing the clique $K_{\ell}$ and minimum degree either 1 or 2, depending of the parity of $g$.
Thus, $\alpha(G_{SR})=2=h+1$.

{\bf Case (c)}:
Suppose that $h\ge2$ and  $g=2k$ is even.
According to Remark \ref{sdimremark2}, $h=2r+s$ and $r+s \le \alpha (G_{SR}) \le  r+s +1$.
Hence, $r=0$, i.e., $C_2(C_g)$ contains no pair of antipodal vertices, if and only if both $h=s$ and 
$s \le \alpha (G_{SR}) \le  s +1$.
Moreover, $\alpha (G_{SR}) =  s $ if and only if $h=k$, which means that 
$\alpha (G_{SR}) =  s +1 =h+1$ if and only if $r=0$ and $h \le k-1$.

{\bf Case (d)}:
Suppose that $h\ge2$ and  $g=2k+1$ is odd.
According to Remark \ref{sdimremark1}, $\left \lceil \frac{h}{2} \right \rceil \le \alpha(G_{SR}) \le h+1$.
In addition, $\alpha (G_{SR}) = h+1$ if and only if both {\colred $G_{SR}[W]=hK_1$ and $C_3(C_g)$ contain} at least one antipodal triple, which also includes the fact that $h \le k-1$.
\end{proof}

\begin{prop}\label{prop.st2}
{\colred
Let $G$ be a proper unicyclic graph of grith $g$ with $\ell$ leaves.
}
If $c_2(C_g)=h$, then 
$\displaystyle \sdim(G)=\ell + \left \lfloor \frac{h}{2} \right \rfloor$
if and only if any of the following conditions hold:

\begin{enumerate}[label=\rm \bf(\arabic*)]

\item  $g$ is even and $h=g-1$.

\item $g$ is odd and $g-2 \le h \le g-1$. 

\end{enumerate}
\label{sdim.ell+h/2}
\end{prop}
\begin{proof} 
According to Theorem \ref{sdim.bounds} and Remarks \ref{sdimremark1} and \ref{sdimremark2}, 
$\sdim(G)=\ell + \left \lfloor \frac{h}{2} \right \rfloor$ if and only if 
$\alpha (G_{SR}) = \left \lceil \frac{h}{2} \right \rceil$.
We distinguish cases:

\vspace{.2cm}
{\bf Case (a)}:
Suppose that $g$ is even.
According to Remark \ref{sdimremark2}, $h=2r+s$ and $r+s \le \alpha (G_{SR}) \le  r+s +1$.
Assume that $h=g-1$. 
Then, $2r+s=2r+2s+2t-1$, which means that $t=0$, $s=1$ and $h=2r+1$.
So, $\alpha (G_{SR}) = r+1= \left \lceil  \frac{h}{2} \right \rceil $.

Conversely, suppose that $\alpha (G_{SR}) = \left \lceil  \frac{h}{2} \right \rceil  = 
r + \left \lceil  \frac{s}{2} \right \rceil $.
Thus, $r+s \le r + \left \lceil  \frac{s}{2} \right \rceil \le  r+s +1$, i.e., 
$s \le  \left \lceil  \frac{s}{2} \right \rceil \le s +1$, which means that 
$s =  \left \lceil  \frac{s}{2} \right \rceil$
since $\left \lceil  \frac{s}{2} \right \rceil < s +1$, for every positive integer $s$.
After noticing that $t=0$ and $g=2r+2s+2t$,
we distinguish cases:

\vspace{.2cm}
{\bf Case (a.1)}:
If $s=0$, then $G=C_g$, a contradiction.

\vspace{.2cm}
{\bf Case (a.2)}:
If $s=1$, then $h=2r+1$ and $g=2r+2$, i.e., $h=g-1$.

\vspace{.2cm}
{\bf Case (b)}:
Suppose that $g$ is odd.
If $g-2 \le h \le g-1$,  then according to Remark \ref{sdimremark1}, in all possible cases,  
$\alpha (G_{SR})= \left \lceil  \frac{h}{2} \right \rceil$.

Conversely, suppose both that $\alpha (G_{SR}) = \left \lceil  \frac{h}{2} \right \rceil $ and $h\le g-3$.

If $\displaystyle G[W]=\sum_{i=1}^{p} P_{h_i}$, then
$$\displaystyle \sum_{i=1}^{p} \alpha(P_{h_i})   \le   \alpha (G_{SR}) \le \sum_{i=1}^{p} \alpha(P_{h_i}) +1.$$
Hence, $\alpha (G_{SR}) = \left \lceil  \frac{h}{2} \right \rceil $ if and only if 
$\displaystyle \alpha (G_{SR}) =  \sum_{i=1}^{p} \alpha(P_{h_i}) = \left \lceil  \frac{h}{2} \right \rceil$.

{\colred 
Notice that, for every pair  $h_i,h_j$ of odd integers, 
$\alpha(P_{h_i})+\alpha(P_{h_j}) = \frac{h_i+h_j}{2} +1$.
Hence,  
$\displaystyle \sum_{i=1}^{p} \alpha(P_{h_i}) = \left \lceil  \frac{h}{2} \right \rceil$ 
if and only if there is at most an index $j \in [p]$ such that $P_j$  is an odd path. We distinguish cases.}

\vspace{.2cm}
{\bf Case (c)}:
There is a unique index $j \in [p]$ such that $P_j$  is an odd path.
We distinguish cases.

\vspace{.2cm}
{\bf Case (c.1)}:
If $p=2$, then we can suppose w.l.o.g. that $P_{h_1}$ is the even path and that between $P_{h_1}$ and $P_{h_2}$ there are in $G_{SR}$ at least two vertices (see Figure \ref{fig:sdim6}, left).
Certainly, it is possible to select a maximum independent set in $P_1$  not containing its last vertex. 
Hence, 
$\displaystyle \alpha (G_{SR}) =  \sum_{i=1}^{p} \alpha(P_{h_i}) + 1 > \left \lceil  \frac{h}{2} \right \rceil$, a contradiction.

\vspace{.2cm}
{\bf Case (c.2)}:
Assume that $p\ge 3$.
Take $P_{h_1}$ and $P_{h_2}$ and suppose w.l.o.g. that they are both even paths and  consecutive in $G_{SR}$ (see Figure \ref{fig:sdim6}, right).
Clearly, it is possible to select a maximum independent set in $P_{h_1}$ (resp., in $P_{h_2}$) not containing its last vertex (resp., its first vertex). 
Hence, 
$\displaystyle \alpha (G_{SR}) =  \sum_{i=1}^{p} \alpha(P_{h_i}) + 1 > \left \lceil  \frac{h}{2} \right \rceil$, a contradiction.

\begin{figure}[!htbp]
\begin{center}
\includegraphics[width=0.78\textwidth]{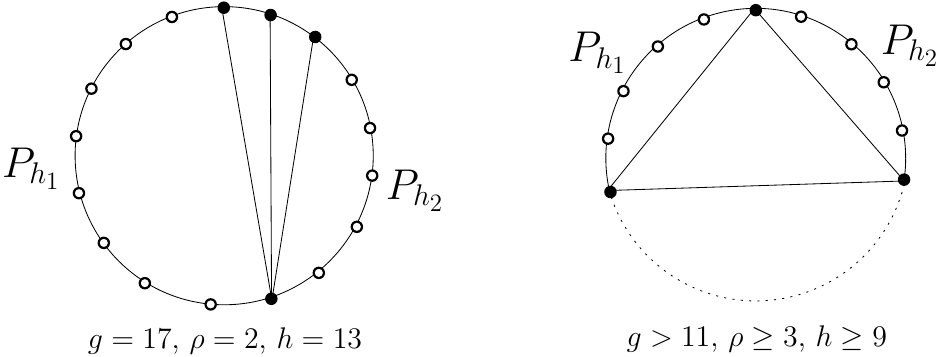}
\end{center}
\caption{Strong locating set of a uniclyclic graph of odd girth $g$ and $h\le g-3$. Left: $p=2$. Right: $p\ge3$. }
\label{fig:sdim6}
\end{figure}

{\bf Case (d)}:
For every index $j \in [p]$,  $P_{h_1}$  is an even path.
Take $P_{h_1}$ and $P_{h_2}$ and suppose w.l.o.g. that they are consecutive in $G_{SR}$ (see Figure \ref{fig:sdim6}, right).
Clearly, it is possible to select a maximum independent set in $P_{h_1}$ (resp., in $P_{h_2}$) not containing its last vertex (resp., its first vertex). 
Hence, 
$\displaystyle \alpha (G_{SR}) =  \sum_{i=1}^{p} \alpha(P_{h_i}) + 1 > \left \lceil  \frac{h}{2} \right \rceil$, a contradiction.
\end{proof}

As for the problem of characterizing the  family of proper unicyclic  graphs $G$ of girth $g$ with at most $\lceil \frac{g}{2} \rceil$  leaves such that $\sdim(G)=\lceil \frac{g}{2} \rceil$, it was solved when $g$ is even in \cite{k20},  while the odd case still remains open.

\begin{prop}
{\colred
Let $G$ be an even  proper unicyclic graph of grith $g$ with $\ell$ leaves.
}
Then, $\sdim(G)=\frac{g}{2}$ if and only if the following conditions hold:

\begin{enumerate}[label=\rm \bf(\arabic*)]

\item  $\rho(G)=0$

\item There is at most one pair of antipodal vertices in $C_g$, each one of degree at least 3.

\end{enumerate}
\label{sdim.g/2}
\end{prop}


Next, we present a result that clearly shows how easy is to compute the strong metric dimension of an even unicyclic graph.

\begin{theorem}\label{th:sd}
{\colred
Let $G$ be a unicyclic graph of even grith $g=2k$ with $\ell$ leaves.
}
If $C_2(C_g)$ contains $r$ antipodal pairs and $C_3(C_g)$ contains $t$ antipodal pairs, then

$$
\sdim(G)=
\left \{ \begin{array}{cl}
\ell + r -1 & {\rm if}\, t \ge 1, \\
\ell + r & {\rm if}\, t=0.  \\
\end{array}\right . 
$$
\label{sdim.even}
\end{theorem}
\begin{proof} 
Let $s$ such that $2s=g-2r-2t$. 
According to Theorem \ref{sdim.partalpha},
 
$$\sdim(G)=|\partial(G)|-\alpha(G_{SR}) = \ell + c_2(C_g) - \alpha(G_{SR}) = \ell +2r + s - \alpha(G_{SR}).$$

As shown in Figure \ref{fig:sdim4} (when $g=8$, $r=2$, $s=1$ and $t=1$), it is easy to check that  
$\alpha(G_{SR})=r+s$ (resp., $\alpha(G_{SR}=r+s+1$) if $t=0$ (resp., $t\ge1$).
\\ 
Hence, $\sdim(G)= \ell +2r + s - r-s = \ell + r$ (resp., $\sdim(G)= \ell +2r + s - r-s -1 = \ell + r -1$)  if $t=0$ (resp., $t\ge1$).
\end{proof}

We finalize this section by showing that the problem of computing the strong metric dimension of an odd unicyclic  graph is  certainly quite more complicated than in the  even case.

\begin{prop}
{\colred
Let $G$ be a unicyclic graph of odd grith $5 \le g=2k+1$ with $\ell$ leaves. 
If $2 \le c_2(C_g)=h \le 2k-2$, then
}

$$
\max\{\ell-1, \ell +h -k\} \le \sdim(G) \le \ell + \Big\lfloor \frac{h}{2}  \Big\rfloor -1 .
$$

Moreover, both bounds are tight.
\label{sdim.odd3}
\end{prop}
\begin{proof} 
The inequality $\ell-1\le \sdim(G)$ and the upper bound are a direct consequence of 
{\colred Lemma \ref{sdimlema2}
}
 and Proposition \ref{sdim.ell+h/2}, respectively.
Hence, to prove the lower bound  it is enough to show that if $ k-1 \le h \le 2k-2$ , then $\alpha(G_{SR}) \le k$, since $\sdim(G) = \ell + h - \alpha(G_{SR})$.

According to Remark \ref{sdimremark1}, suppose that $G_{SR}[W]$ is the disjoint union of $p$ paths, in such a way that $p=\mu+\eta$, there are $\mu$ odd paths and $\eta$ even paths.
If $\displaystyle G_{SR}[W] = \sum_{\kappa=1}^{\mu} P_{i_{\kappa}} + \sum_{\kappa=1}^{\eta} P_{j_{\kappa}}$, then  
$$\displaystyle \sum_{\kappa=1}^{\mu} \alpha(P_{i_{\kappa}}) + \sum_{\kappa=1}^{\eta} \alpha(P_{j_{\kappa}})
\le
\alpha(G_{SR}) 
\le 
\sum_{\kappa=1}^{\mu} \alpha(P_{i_{\kappa}}) + \sum_{\kappa=1}^{\eta} \alpha(P_{j_{\kappa}}) +1.$$
If for every $\kappa \in [\mu]$ (resp., $\kappa \in [\eta]$), 
$i_{\kappa}=2r_{\kappa}+1$ (resp., $j_{\kappa}=2s_{\kappa}$ ), then 

$$\displaystyle \sum_{\kappa=1}^{\mu} r_{\kappa} + \sum_{\kappa=1}^{\eta} s_{\kappa} + \mu
\le
\alpha(G_{SR}) 
\le 
\sum_{\kappa=1}^{\mu} r_{\kappa} + \sum_{\kappa=1}^{\eta} s_{\kappa} + \mu +1.$$

and

$$
\displaystyle
g
\ge
\sum_{\kappa=1}^{\mu} 2r_{\kappa} + \mu + \sum_{\kappa=1}^{\eta} 2s_{\kappa} + p
= 2\left[ 
\sum_{\kappa=1}^{\mu} r_{\kappa} + \sum_{\kappa=1}^{\eta} 2s_{\kappa} + \mu 
\right] + \eta.
$$

We distinguish cases.

\vspace{.2cm}
{\bf Case (1)}:
$\eta \ge 2$.
In this case,  
$k \ge \sum_{\kappa=1}^{\mu} r_{\kappa} + \sum_{\kappa=1}^{\eta} 2s_{\kappa} + \mu +1  \ge \alpha(G_{SR})$, since $g$ is odd.

\vspace{.2cm}
{\bf Case (2)}:
$g-h \ge p +2$.
In this case,
$$\displaystyle g \ge h + p +2 = 
h + \mu + \eta +2 =
2\left[ 
\sum_{\kappa=1}^{\mu} r_{\kappa} + \sum_{\kappa=1}^{\eta} s_{\kappa} + \mu +1
\right] + \eta.$$
Hence,
$k \ge \sum_{\kappa=1}^{\mu} r_{\kappa} + \sum_{\kappa=1}^{\eta} s_{\kappa} + \mu +1 \ge \alpha(G_{SR})$, since $g$ is odd.

\vspace{.2cm}
{\bf Case (3)}:
$p \le g-h \le p +1$ and $0 \le \eta \le 1$.
In this case, 
$\displaystyle \alpha(G_{SR})= \sum_{\kappa=1}^{\mu} r_{\kappa} + \sum_{\kappa=1}^{\eta} s_{\kappa} +\mu$.

Hence, $k \ge \sum_{\kappa=1}^{\mu} r_{\kappa} + \sum_{\kappa=1}^{\eta} s_{\kappa} + \mu  = \alpha(G_{SR})$, since $g$ is odd.
%
%
%

To show the tightness of the upper bound, consider the case  $h=g-3=2k-2$ and $p=\eta=1$, and notice that
{\colred $\alpha(G_{SR})=\alpha(P_{2k-2})+1=k$.}
Hence, $\sdim(G)=\ell + h -\alpha(G_{SR}) = \ell+  2k-2 - k = \ell + k -2 = \ell + \Big\lfloor \frac{h}{2}  \Big\rfloor -1$.

The tightness of the lower bound when  $0 \le h \le k-1$ is a direct consequence of Proposition \ref{sdim.ell-1}.

If $h=k$, consider the case $\displaystyle G_{SR}[W] = kK_{1}$ and notice that $\alpha(G_{SR})=k$.
If $k+1 \le h \le 2k-2$, then consider the case
$\displaystyle G_{SR}[W] = (g-h-1)K_{1} + P_{2h-2k}$.
Then, $\alpha(G_{SR})= g-h-1 + \alpha(P_{2h-2k}) = 2k+1-h-1+h-k=k$.
Thus, if $k\le h \le 2k-2$, then $\sdim(G)=\ell + h -\alpha(G_{SR})= \ell + h -k$.
\end{proof}

\section{Dominating metric dimension}\label{ddim}

After the introduction of metric dimension, it seems natural to study when a metric-locating set verifies additional properties, and among all the vertex-set properties, one of the most important is domination. 
Thus, R. C. Brigham et al. in \cite{bcdz03}, and also and independently M. A. Henning and O. R.  Oellermann in \cite{ho04}, defined the following parameter.

\begin{defi}
{\rm A set of vertices of a graph $G$ is called \emph{metric-locating-dominating} (an \emph{MLD-set} for short) if it is both locating and dominating.}
\end{defi}

An MLD-set of  minimum cardinality is called  an \emph{MLD-basis} of $G$. 
The \emph{dominating metric dimension} of $G$, denoted by $\ddim(G)$,  is the  cardinality of an MLD-basis. 
To know more about this parameter see mainly  \cite{bcdz03,ho04,chmpp13} and also \cite{ghm18,hmp14,zab22}.

As usual, the dominating metric dimension of trees is well-known.

\begin{theorem} 
{\rm\cite{ho04}}           
Let $T$ be a tree with $\ell(T)$ leaves and $s(T)$ support vertices.\\
Then,  $\ddim(T)=\gamma(T)+ \ell(T) - s(T)$.
\label{thmho04.1}
\end{theorem}

\begin{cor} 
{\rm\cite{ho04}}           
Let $T$ be a tree.
Then,  $\ddim(T)=\gamma(T)$ if and only if $T$ contains no strong support vertex.
\label{corho04.1}
\end{cor}

In particular, $\ddim(P_n)=\gamma(P_n)=\lceil \frac{n}{3}  \rceil$. 
The same reasonIng as in trees can be applied in general graphs to obtain a lower bound.

\begin{prop}    
Let $G$ be a graph with $\ell(G)$ leaves and $s(G)$ support vertices.\\
Then,  $\ddim(G) \ge \gamma(G)+ \ell(G) - s(G)$.
\label{lmddim3}
\end{prop}
\begin{proof}
Consider the set $S(G)$ of  support vertices of $G$.
{\colred
For every vertex $v\in S(G)$, let $L_v$ denote the set of leaves adjacent to $v$, i.e.,
$L_v=N(v)\cap {\cal L}(G)$.
}
Let $D$ be a minimum MLD-set of $G$, containing as few leaves as possible.
This means  that, for every vertex $v \in S(G)$, $D$ contains all except one leaf $v'$ adjacent to $v$, as well as vertex $v$.
So, the set $\displaystyle D'= D -\Big( \cup_{v\in S(G)} (L_v-v')  \Big)$
is a dominating set of $G$.
Hence,
$$\gamma(G) \le |D'| \le  |D| - \sum_{v\in S(G)} (|L_v|-1) = \ddim(G) -\ell(G) +s(G),$$
and the inequality $\ddim(G) \ge \gamma(G)+ \ell(G) - s(G)$ follows.
\end{proof}

As an immediate consequence of this result the following one holds.

\begin{lemma}    
For every graph  $G$, $\gamma(G) \le \ddim(G)$.
Moreover, if  $\ddim(G) = \gamma(G)$, then $G$ has no strong support vertex.
\label{lmddim2}
\end{lemma}

As a consequence, there exist several situations in which a dominating set is also an MLD-set.

\begin{figure}[ht]
\begin{center}
\includegraphics[width=0.85\textwidth]{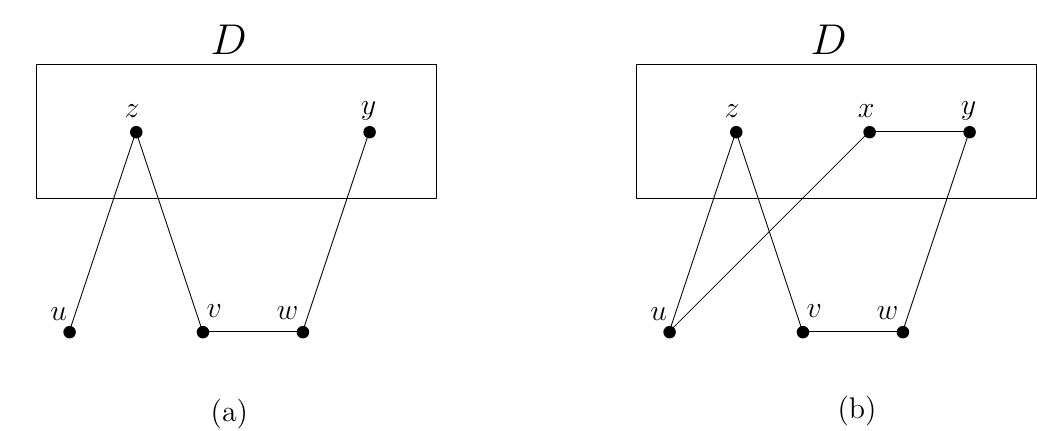}
\caption{$D$ is a dominating set of a unicyclic graph with no strong support vertex.}
\label{uvzwyx}
\end{center}
\end{figure}

\begin{prop}  
{\colred
Let $G$ be a unicyclic graph of girth $g$ with no strong support vertex.
}
If $D$ is a dominating set of $G$, then $D$ is an MLD-set of $G$ whenever any of the following conditions hold.

\begin{enumerate}[label=\rm \bf(\arabic*)]

\item $g \not\in \{3,4,6\}$.

\item $g=3$ and $|D\cap V(C_3)|\ge2$.

\item $g \in \{4,6\}$ and $|D\cap V(C_g)|\ge3$.

\end{enumerate}

\label{lmddim1}
\end{prop}
\begin{proof}
On the contrary, let $D$ be a dominating set $D$ verifying the conditions of the hypothesis, and suppose $D$ is not an MLD-set, i.e, there exist $u,v \in V(G)-D$ such that $N(u)\cap D = N(v) \cap D$.
Notice that there is a unique vertex $z$ such that $N(u)\cap D = N(v) \cap D =\{z\}$, 
since either $g\neq 4$ or $g=4$ and $|D\cap V(C_4)|\ge3$.
As $z$ is not a strong support vertex, then we can suppose w.l.o.g. that $deg(v)\ge 2$.
Take a vertex $w \in N(v)-z$.
Notice that $w \not\in D$, as $N(v) \cap D =\{z\}$.
Let $y\in N(w) \cap D$.
Observe that $y\neq z$, as either $g\neq 3$ or $g=3$ and $|D\cap V(C_3)|\ge2$.
Certainly, $d(v,y)=2$ and it is also clear that $d(u,y)>2$, unless $g=6$ (see Figure \ref{uvzwyx} (a)).
If $g=6$ and $d(u,y)=2$, take the  vertex $x\in N(u)\cap N(y)$ and notice first that necessarily it belongs to $D$ since $|D\cap V(C_6)|\ge3$, and second that $d(u,x)=1$ and $d(v,x)=3$ (see Figure \ref{uvzwyx} (b)).
%
\end{proof}

From the preceding two results, the following one  immediately follows.

\begin{cor}    
Let $G$ be a unicyclic graph of girth $g$.
If $g \not\in \{3,4,6\}$,
then  $\ddim(G)=\gamma(G)$ if and only if $G$ contains no strong support vertex.
\label{corddim1}
\end{cor}

In particular, if $g \not\in \{3,4,6\}$, then  $\ddim(C_g)= \gamma(C_g)=\lceil \frac{g}{3} \rceil$. 
The previous results lead us to obtain the same characterization as trees for unicyclic graphs except for the cases $g\in\{3,4,6\}$.

\begin{theorem}    
Let $G$ be a unicyclic graph of girth $g$ with $\ell(G)$ leaves and $s(G)$ support vertices.
If $g \not\in \{3,4,6\}$, then  $\ddim(G)=\gamma(G)+ \ell(G) - s(G)$.
\label{thmddim1}
\end{theorem}
\begin{proof}
According to Proposition \ref{lmddim3}, $\ddim(G) \ge \gamma(G)+ \ell(G) - s(G)$.
To prove the other inequality,
consider the set  $S(G)$ of strong support vertices of $G$.
If $|S(G)|$=0, then $\ell(G) = s(G)$, and according to Corollary \ref{corddim1}, the equality holds.
Suppose thus that $|S(G)|\ge1$.
For each vertex $v \in S(G)$, take a vertex $v' \in L_v = N(v)\cap {\cal L}(G)$.
Consider the unicyclic graph $G'$ obtained from $G$ by deleting, for every vertex $v \in S(G)$, the set $L_v -v'$.
Hence, $G'$ has no strong support vertex.
Let $D'$ be a minimum dominating set of $G'$ containing all the support vertices of $G'$.
According to Proposition \ref{lmddim1}, $D'$ is an MLD-set of $G'$.
Thus, 
$\displaystyle  D' \cup \Big( \cup_{v\in S(G)} (L_v-v') \Big)$ 
is an MLD-set of $G$, and so
$$\displaystyle \ddim(G) \le |D'| + \sum_{v\in S(G)} (|L_v|-1) = \gamma(G) +\ell(G) - s(G).$$
\end{proof}

The cases $g=3$, $g=4$ and $g=6$ are approached in the following two results.

\begin{figure}[ht]
\begin{center}
\includegraphics[width=0.80\textwidth]{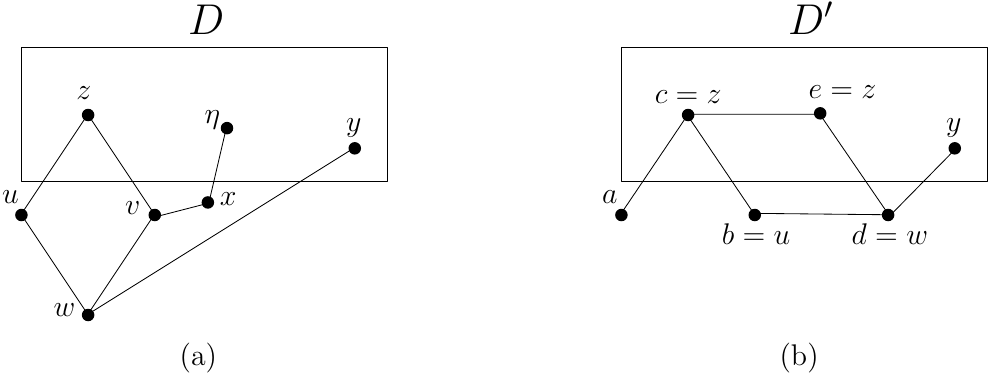}
\caption{$D$ is a dominating set of a unicyclic graph of grith $g=4$, with no strong support vertex.}
\label{abcdeyeta}
\end{center}
\end{figure}

\begin{lemma}
{  \colred
Let $G$ be a unicyclic graph of girth $g=4$ with no strong support vertex.
}
If $D$ is a minimum dominating set of $G$ and $|D\cap V(C_g)|=1$, then $\gamma(G) \le \ddim(G) \le \gamma(G)+1$.
\label{lmddim4}
\end{lemma}
\begin{proof}
Let $u,v \in V(G)-D$.
Suppose that $N(u)\cap D = N(v) \cap D$.
Notice that there is a unique vertex $z$ such that $N(u)\cap D = N(v) \cap D =\{z\}$,  since $|D\cap V(C_4)|=1$.
As $z$ is not a strong support vertex, then we can suppose w.l.o.g. that $deg(v)\ge 2$.
Take a vertex $w \in N(v)-z$ and notice that $w \not\in D$, as $N(v) \cap D =\{z\}$.
Let $y\in N(w) \cap D$.
Observe that $y\neq z$, as $g\neq 3$.
Certainly, $d(v,y)=2$ 
Suppose that $d(u,y)=2$, as otherwise we are done.
This means that $w\in N(u)$.
We distinguish cases.

\vspace{.2cm}
{\bf Case (a)}:
If $deg(v)\ge3$, take $x \in N(v)\setminus \{w,z\}$.
Let $\eta \in D \cap N(x)$ and notice that $d(v,\eta)=2$ and $d(u,\eta)=4$ (see Figure \ref{abcdeyeta} (a)).

\vspace{2cm}
{\bf Case (b)}:
If $deg(u)=deg(v)=2$, consider the set $D'=D+v$.
In order to prove that $D'$  is an MLD-set of $G$, take a pair of vertices $a,b \in V(G)-D'$.
If $N(a)\cap D = N(b) \cap D$, then there is a unique vertex $c$ such that $N(a)\cap D = N(b) \cap D =\{c\}$.
As $c$ is not a strong support vertex, then we can suppose w.l.o.g. that $deg(b)\ge 2$.
Take a vertex $d \in N(b)-c$ and notice that $d\not\in D$, as $N(b) \cap D =\{c\}$.
Let $e\in N(d) \cap D$ and check that $e\neq e$.
Certainly, $d(b,e)=2$.
Suppose that $d(a,e)=2$, as otherwise we are done.
This means that $e\in N(c)$ and also that $c=z$, $e=v$, $b=u$ and $d=w$.
Hence, $d(b,y)=d(u,y)=2$ and $d(a,y)=4$ (see Figure \ref{abcdeyeta} (b)).
\end{proof}

By considering the previous Lemma along with  items (2) and (3) of Proposition \ref{lmddim1}, and reasoning in  a similar way as in the proof Theorem \ref{thmddim1}, the following result holds.

\begin{prop}    
Let $G$ be a unicyclic graph of girth $g$ with $\ell(T)$ leaves and $s(T)$ support vertices.
If $g \in \{3,4,6\}$, then  $\gamma(G)+ \ell(T) - s(T) \le \ddim(G) \le \gamma(G)+ \ell(T) - s(T)+1$.
\label{thmddim2}
\end{prop}

\begin{figure}[ht]
\begin{center}
\includegraphics[width=0.97\textwidth]{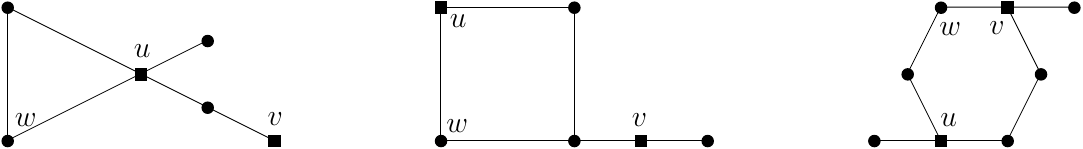}
\caption{ In all cases, $\{u,v\}$ is  a dominating set, $\gamma(G)=2$, $\{u,v,w\}$ is  an MLD-basis and
$\ddim(G)=3$.}
\label{g346}
\end{center}
\end{figure}

As for tightness of the  above upper bound, in some cases it is needed to add  an additional vertex to have an MLD-set starting from a minimum dominating set and the set of leaves hanging on the strong support vertices (see Figure \ref{g346}, for some examples).

\section{Other metric locating parameters}\label{omlp}

In this section, we revise the situation of other metric locating parameters regarding unicyclic graphs. 
For the last two of them, a completed characterization has been obtained by other authors. 
However, for the three first parameters that characterization is far from being obtained, and here we described the advances in that matter.

\subsection{Fault-tolerant metric dimension}\label{dim2}

A  set of vertices $S$ of a graph  $G$ is called \emph{$2$-locating} if  every pair of distinct  vertices  $x,y \in V (G)$ is resolved by at least $2$ elements of $S$.

A $2$-locating set of  minimum cardinality is called  a \emph{$2$-metric basis} of $G$. 
The \emph{$2$-metric dimension} of $G$, denoted by $\dim_2(G)$,  is the  cardinality of a $2$-metric basis.

This parameter was formally introduced by C. Hernando et al. in \cite{hmsw08}. 
In this work, 2-metric locating sets were called \emph{fault-tolerant sets}, and the 2-metric dimension was  named the \emph{fault-tolerant metric dimension} of $G$.
To know more about this parameter see also  \cite{cjs10,jscs09}.

It is a routine exercise to check that  $\dim_2(P_n)=2$ and $\dim_2(C_n)=3$.

\begin{theorem}
{\rm\cite{hmsw08}}
Let  $T$ be  a tree  having $\ell_s$ strong leaves. 
Then, $\dim_2(T)  = \ell_s$.
\end{theorem}

\subsection{$k$-metric dimension}\label{dimk}

Let $k,n$ be a pair of vertices such that $2 \le k \le n-1$.
A  set of vertices $S$ of a graph  $G$ of order $n$ is called \emph{$k$-locating} if  every pair of distinct  vertices  $x,y \in V (G)$ is resolved by at least $k$ elements of $S$.

A $k$-locating set of  minimum cardinality is called  a \emph{$k$-metric basis} of $G$. 
The \emph{$k$-metric dimension} of $G$, denoted by $\dim_k(G)$,  is the  cardinality of a $k$-metric basis. 

This parameter was formally introduced by A. Estrada-Moreno in \cite{ery15}. 
The case $k=2$ was previously introduced by C. Hernando et al. in \cite{hmsw08} (see  previous subsection). 
To know more about this parameter see also \cite{eyr14,eyr16,eyr16.2,mst22}.

A  graph  $G$ is called 
\emph{$k$-metric dimensional}, if $k$ is the largest integer such
that $G$ contains a $k$-locating set.

\begin{prop}
{\rm \cite{e21, ery15}}
Let  $k,n\ge3$ be integers such that $k+1 \le n$.
Then, 
\begin{itemize}

\item $P_n$ is $(n-1)$-dimensional.

\item $\dim_k(P_n)=k+1$.

\item If $n$ is odd (resp., even) , then $C_n$ is $(n-1)$-dimensional (resp., $(n-2)$-dimensional).

\item $\dim_k(C_n)=
\left \{ 
\begin{array}{cl}
k+2 & {\rm if}\, n  \, {\rm is}\, {\rm even}\, {\rm and}\, \frac{n}{2} \le k \le n-2, \\
k+1 & {\rm otherwise}.  \\
\end{array}
\right .$.

\end{itemize}

\label{propc4.f1}
\end{prop}

\begin{theorem}
{\rm \cite{ery15}}
{\colred Let  $G$ be a $k$-metric dimensional graph of order $n\ge3$.}
Then,

\begin{itemize}

\item $2\leq k\leq n-1$.

\item $k=2$ if and only if $G$ has  twin vertices.

\end{itemize}

\end{theorem}

Let ${\cal M}(T)$ denote the set of strong exterior major vertices of a tree $T$.
For every vertex $w \in  {\cal M}(T)$,
let $Ter(w)$ and $ter(w)$ denote the set and the number of terminal vertices of $w$, respectively.

If $Ter(w)=\{u_1,\ldots,u_h\}$, then 
$\displaystyle l(w)=\min_{i\in[h]}d(u_i,w)$
and
$\displaystyle \zeta(w)=\min_{\stackrel{i,j\in[h]}{i\neq  j}}[d(u_i,w)+d(w,w_j)]$. Finally, $\displaystyle  \zeta(T)=\min_{w \in {\cal M}(T)} \zeta(w)$ and

$I_r(w)=
\left \{ 
\begin{array}{cl}

[ter(w)-1][r-l(w)]+l(w) & {\rm if}\,   l(w) \le \, \lfloor \frac{r}{2} \rfloor\\

[ter(w)-1]\lceil \frac{r}{2} \rceil +  \lfloor \frac{r}{2} \rfloor & {\rm otherwise}.  \\

\end{array}
\right .
$.

\begin{theorem}
{\rm \cite{ery15}}
Let  $T$ be a $k$-metric dimensional tree, other than a path.
Then, 
\begin{itemize}

\item $k=\zeta(T)$.

\item For every $r \in [\zeta(T)]$,
$\dim_r (T) = \sum_{w \in {{\cal M}(T)}} I_r(w) $.

\end{itemize}

\end{theorem}

In \cite{e21}, some partial results for unicyclic graphs are presented. 
For example, it was given a closed formula for the $k$-metric dimension of a unicyclic graph $G$ of girth $g$, whenever $\rho(G)=0$ and $c_2(C_g)=g-1$, where $C_g$ is the cycle of $G$.

\subsection{Edge metric dimension}\label{edim}

Let $G=(V,E)$ be a graph.
A vertex $v\in V$ is said to resolve two edges $e,f \in E$ if $d_G(v,e)\neq d_G(v,f)$.
A  set of vertices $S$ of a graph $G$ is called \emph{edge locating} if  every pair of distinct  edges  $e,f\in E$ is resolved by a vertex of $S$.

An edge locating set of  minimum cardinality is called  a \emph{edge metric basis} of $G$. 
The \emph{edge metric dimension} of $G$, denoted by $\edim(G)$,  is the  cardinality of an edge metric basis. 
This parameter was formally introduced by A. Kelenc et al. in \cite{kty18}. 
To know more about this parameter see also \cite{azszas22,g20,kmmsy21,py20,ss22.2,ss22.1}.

It is a routine exercise to check that $\edim(P_n)=1$ and $\edim(C_g)=2$.

\begin{theorem}
{\rm\cite{kty18}}
{ \colred
For every tree  $T$, $\edim(T)  = \ell(T)-\lambda(T)$.
}
\end{theorem}

\begin{theorem}
{\rm \cite{ss22}}
{\colred
Let  $G$ be a unicyclic graph  such that $\ell(G)=\ell$ leaves, $\lambda(G)=\lambda$  and $\rho(G)=\rho$-
Then,
}  $$\ell - \lambda + \max\{2-\rho,0\} \le \edim(G ) \le  \ell - \lambda + \max\{2-\rho,0\} +1.$$
\end{theorem}

\begin{theorem}
{\rm \cite{wyc22}}
Let $G$ be a unicyclic graph.
If  $e\in E(G)$ and $G-e$ is connected, then $\edim(G) \le \edim(G-e) +1$.
\end{theorem}

\begin{theorem}
{\rm \cite{ss22}}
{\colred
If  $G$ is a unicyclic graph, then $|\dim(G)-\edim(G)| \le 1$.
}
\end{theorem}

Moreover, a  characterization of the sets of unicyclic graphs  attaining each of the three possible equalities  was also given in \cite{ss22}.
In particular, it was proved that if $G$ is an odd (resp., even) unicyclic graph, then $\dim(G)\le\edim(G)$ (resp.,  $\dim(G)\ge\edim(G)$).

\subsection{Mixed metric dimension}\label{mdim}

A  set of vertices $S$ of a graph $G$ is called \emph{mixed locating} if  every pair of distinct  elements (vertices or edges)   is resolved by a vertex  of $S$.

A mixed locating set of  minimum cardinality is called  a \emph{mixed metric basis} of $G$. 
The \emph{mixed metric dimension} of $G$, denoted by $\mdim(G)$,  is the  cardinality of mixed  metric basis. 
This parameter was formally introduced by A. Kelenc et al. in \cite{kkty17}. 
To know more about this parameter see  \cite{kkty17,ss21,ss21.2}.

It is straightforward to check that, for every integer $n\ge3$,  $\mdim(P_n)=2$ and $\mdim(C_n)=3$.

\begin{theorem}{\rm \cite{kkty17}}
Let $G$ be a graph of order $n\ge2$.
Then, $\mdim(G)=2$ if and only if $G$ is the path $P_n$.
\end{theorem}

\begin{theorem}{\rm \cite{kkty17}}
For any tree $T$, $\mdim(T)=\ell(T)$.
\end{theorem}

\begin{theorem}{\rm \cite{ss21.2}}
Let $G$ be a unicyclic graph whose cycle $C_g$ has $t$ root vertices.
Then,
$$\mdim(G)=\ell(G)+\max\{3-t,0\}+\epsilon$$
where $\epsilon=1$ if $t\ge3$ and  there is no geodesic triple of root vertices on $C_g$, while $\epsilon=0$ otherwise.
\label{thm1mdim}
\end{theorem}

\subsection{Local metric dimension}\label{ldim}

\begin{def} \label{dls} 
A  set of vertices $S$ of a graph $G$ is called \emph{local locating} if,  for every pair of adjacent vertices $x,y \in V(G)$, there is a vertex $v \in S$, such that 
{\colred $d_G(x,v) \neq d_G(y,V)$.}
\end{def}

In other words,
$S$ is a local locating set of $G$ if, for every pair of adjacent  vertices  $x,y \in E (G)$, $r(x|S) \neq r(y|S)$.

A local locating set of  minimum cardinality is called  a \emph{local metric basis} of $G$. 
The \emph{local metric dimension} of $G$, denoted by $\ldim(G)$,  is the  cardinality of a local metric basis. 
This parameter was formally introduced by F. Okamoto et al. in \cite{opz10}.
To know more about this parameter see also  \cite{berr19,rbg16,rgb15}.

\begin{theorem}
{\rm \cite{opz10}}
{\colred
 Let $G$ be a graph.
 Then, $\ldim(G) =1 $ if and only if it is bipartite.
 }
 \label{th:ldimbipartite}
\end{theorem}

In particular, the local metric dimension of paths, trees and even unicyclic graphs is one.

\begin{theorem}{\rm \cite{frsb22}}
{\colred
Let $G$ be a unicyclic graph. Then,
}
$$ \ldim(G)=\left\{\begin{array}{cr}
1 & g\textrm{ is even,}\\
2 & g\textrm{ is odd.}
\end{array}\right.$$
\label{thmldim1}
\end{theorem}
%

\section{Conclusions and Further work}\label{cfw}

The primary purpose of this paper has been to survey the {\colred state-of-the-art } of the main  metric locating parameters for the family of pseudotrees, i.e., for paths, cycles, trees and unicyclic graphs.

During the process of writing this survey, the authors noticed that regarding paths, cycles and trees, everything is known about the nine parameters analysed here (see Table \ref{megatable1}). 
However, a complete characterization has not been obtained for unicyclic graphs in all the metric locations that we have surveyed, so we set out to contribute with a number of new results.
As a consequence, we have been able to make significant progress in four of the nine parameters considered in this work.

\begin{table}[h]
  \begin{center}
   
$
\begin{array}{c|cccc|ccccc}

\hline

G & {\dmd}  & {\dim} & {\sdim} & {\ddim} & {\dim_2} & {\dim_k} & {\edim} & {\mdim}   & {\ldim} 
\\

\hline 

P_n & 2 &  1  &  1  &  \lceil \frac{n}{3}  \rceil &  2  &  k+1   &   1  &  2   & 1 
\\ 
 &&&&&  &&&& \\

C_n & 2,3 &  2 &  \lceil \frac{n}{2}  \rceil & \lceil \frac{n}{3}  \rceil    &  3  &  \begin{array}{c} k+1,\\k+2 \end{array}   &   2  &  3  &  1,2  
\\ 
&&&&&  &&&& \\

T_n  & \ell &  \ell-\lambda & \ell-1 &  \gamma+\ell-s  &  \ell_s &  \cite{ery15}   &   \ell-\lambda  &  \ell  &  1   
\\ 

\hline

\end{array}
$

\end{center}
\caption{ \label{megatable1} Metric locating parameters of paths, cycles and trees of order $n\ge4$.}
\end{table}

To be more precise, we have completed the work for the doubly metric dimension and the dominating metric dimension, and we also have contributed with significant new results for the metric dimension and the strong metric dimension.

As for the remaining five parameters, for two of them (the mixed metric dimension and the local metric dimension) the values {\colred have already been computed}, while for the remaining three parameters (the edge metric dimension, the fault-tolerant metric dimension and the  $k$-metric dimension), there are plenty of work to be done to get both bounds and exact values.

For further details, see Table \ref{megatable2}  and Table \ref{megatable3}, 
having in mind that  $G$ denotes a proper unicyclic graph  with girth $g$,  $\ell$ leaves, $\ell_s$ strong leaves, $s$ support vertices, $\lambda$ exterior major vertices, $\rho$ branch-active vertices, $c_2$ trivial vertices and $t$ root vertices. 
Moreover, $\hat{\rho}=max\{2-\rho,0\}$, $\hat{t}=max\{3-t,0\}$ and $\mu(g,\ell)=\max\{\lceil \frac{g}{2} \rceil, \ell-1\} $.

\begin{table}[h]
  \begin{center}
    
$
\begin{array}{c|cc|cc}

\hline

{\rm  \, Bounds} & {\dmd}  &   {\ddim} &  {\mdim}   & {\ldim} 
\\

\hline 

{\rm Lower} & \ell &    \gamma + \ell - s &  \ell + \hat{t}   & 1 
\\ 
 &&&  & \\

{\rm Upper} & \ell +2 &    \gamma + \ell - s +1   &   \ell + \hat{t}+1  &  2  
\\ 
&&&  & \\

\hline 

{\rm Exact \, values}  & {\rm Theorem} \, \ref{dmd.unic1} &   {\rm Theorem} \, \ref{thmddim1}  &   {\rm Theorem} \,  \ref{thm1mdim}  &  {\rm Theorem} \, \ref{thmldim1}   
\\ 

\hline

\end{array}
$

 \end{center}
  \caption{\label{megatable2} Metric locating parameters of proper unicyclic graphs (no open problems).}
\end{table}

\begin{table}[h]
  \begin{center}
$
\begin{array}{c|cc|ccc}

\hline

{\rm  \, Bounds} &  {\dim} & {\sdim} &  {\dim_2}  &  {\dim_k} &  {\edim} 
\\

\hline 

{\rm Lower} &  \ell -\lambda +\hat{\rho} &  \mu(g,\ell)  &  -- &   --
&        \ell -\lambda +\hat{\rho}  
\\ 
 &&  &&& \\

{\rm Upper} &  \ell -\lambda +\hat{\rho}+1 &  \ell + \lfloor \frac{c_2}{2}  \rfloor & --  &--
&  \ell -\lambda +\hat{\rho}+1  
\\ 
&&  &&& \\

\hline

{\rm Exact \, values}  &   \begin{array}{c} Cor. \ref{cor:xulo.odd}, \ref{cor:xulo.even}\\and~\cite{ss22}\end{array} & \begin{array}{c}Prop.~\ref{prop.st1},\ref{prop.st2}\\ and~ Th.~\ref{th:sd}\end{array}  & -- & --  &    \cite{ss22}  
\\ 

\hline

\end{array}
$

 \end{center}
  \caption{\label{megatable3} Metric locating parameters of proper unicyclic graphs (with open problems).}
\end{table}

We conclude with a list  of  open problems.

{\rm \bf  Open Problem 1}:
Starting mainly both  from \cite{ss22} and from Section \ref{dim} of this paper, characterizing the family of proper odd unicyclic graphs $G$ such that $\dim(G)=\ell - \lambda + \max\{2-\rho,0\} $.

{\rm \bf  Open Problem 2}:
Starting mainly both  from \cite{ss22} and from Section \ref{dim} of this paper, characterizing the family of proper even unicyclic graphs $G$ such that $\dim(G)=\ell - \lambda + \max\{2-\rho,0\}$.

{\rm \bf  Open Problem 3}:
Starting mainly both  from \cite{k20} and from Section \ref{sdim} of this paper, characterizing the  family of proper odd unicyclic  graphs $G$ of girth $g$ with at most $\lceil \frac{g}{2} \rceil$  leaves such that $\sdim(G)=\lceil \frac{g}{2} \rceil$.

{\rm \bf  Open Problem 4}:
Starting mainly both  from \cite{k20} and from Section \ref{sdim} of this paper, either obtaining a result 
similar to the one displayed in Theorem \ref{sdim.even} when $g$ is odd, or at least characterizing the families of graphs achieving the bounds shown in Proposition \ref{sdim.odd3}.

{\rm \bf  Open Problem 5}:
Starting mainly both  from \cite{ss22} and from Section \ref{dim} of this paper,
characterizing the family of proper unicyclic graphs $G$ such that $\edim(G)=\ell - \lambda + \max\{2-\rho,0\} $, or at least obtaining similar results to those appearing in subsections \ref{dim.odduni} and \ref{dim.evenuni}.

{\rm \bf  Open Problem 6}:
Starting mainly from \cite{e21}, \cite{ery15} and \cite{hmsw08}, obtain tight, both lower and upper bounds, of the fault-tolerant metric dimension of proper unicyclic graphs. 

{\rm \bf  Open Problem 7}:
Starting mainly from \cite{e21} and \cite{ery15}, obtain tight, both lower and upper bounds, of the $k$-metric dimension of proper unicyclic graphs without strong support vertices, for every integer $k\ge3$. 




\end{document}